\documentclass[12pt]{article}
\usepackage{amstext,amsfonts,amssymb,latexsym,graphics}
\newtheorem{thm}[equation]{Theorem}
\newtheorem{lemma}[equation]{Lemma}
\newtheorem{prop}[equation]{Proposition}

\newenvironment{pf}
{\begin{trivlist} \item \noindent{\sc Proof. }}
{{\hfill $\Box$} \end{trivlist}}
\newenvironment{pfof}[1]
{\begin{trivlist} \item \noindent{\sc Proof of\ #1. }}
{{\hfill $\Box$} \end{trivlist}}

\newcommand{\zt}{{\mathbb Z}^2}

\newcommand{\rtp}{{\mathbb R}^2_+}
\newcommand{\vc}[1]{\mathbf{#1}}
\newcommand{\lng}{\,\mbox{long}}
\newcommand{\sht}{\,\mbox{short}}
\newcommand{\sem}{\phi}
\newcommand{\nnorm}{N_{\text{\rm norm}}}
\newcommand{\nsplit}{N_{\text{\rm split}}}
\newcommand{\intint}{\int\hspace{-0.7em}\int}

\newcounter{mycount}
\newenvironment{mylist}{\begin{list}{{\rm (\roman{mycount})}}%
{\usecounter{mycount}\itemsep 0pt}}{\end{list}}

\newcommand{\dof}{\em}

\begin{document}

\title{Sharp Metastability Threshold for Two-Dimensional Bootstrap Percolation}
\author{Alexander E. Holroyd}
\date{May 8, 2002 (updated \today)}

\maketitle

\begin{abstract}
In the bootstrap percolation model, sites in an $L$ by $L$ square are initially independently declared active with probability $p$.  At each time step, an inactive site becomes active if at least two of its four neighbours are active.  We study the behaviour as $p\rightarrow 0$ and $L\rightarrow \infty$ simultaneously of the probability $I(L,p)$ that the entire square is eventually active.  We prove that $I(L,p)\rightarrow 1$ if $\liminf \;p\log L>\lambda$, and $I(L,p)\rightarrow 0$ if $\limsup \;p\log L<\lambda$, where $\lambda=\pi^2/18$.
We prove the same behaviour, with the same threshold $\lambda$, for the probability $J(L,p)$ that a site is active by time $L$ in the process on the infinite lattice.  The same results hold for the so-called modified bootstrap percolation model, but with threshold $\lambda'=\pi^2/6$.  The existence of the thresholds $\lambda,\lambda'$ settles a conjecture of Aizenman and Lebowitz \cite{aizenman-leb}, while the determination of their values corrects numerical predictions of Adler, Stauffer and Aharony \cite{adler2}. 
\renewcommand{\thefootnote}{}
\noindent\footnote{\hspace{-2em}{\bf\noindent Key words:} bootstrap percolation, cellular automaton, metastability, finite-size scaling}
\noindent\footnote{\hspace{-2em}{\bf\noindent 2000 Mathematics Subject
Classifications:} Primary 60K35; Secondary 82B43}
\noindent\footnote{\hspace{-2em}{\bf Address:} UCLA Department of Mathematics, CA 90095-1555, USA}
\noindent\footnote{\hspace{-2em}Research funded in part by NSF Grant DMS--0072398}
\renewcommand{\thefootnote}{\arabic{footnote}}
\end{abstract}

\section{Introduction}
\label{intro}

We consider the bootstrap percolation process in two dimensions.
Let $\zt=\{x=(x_1,x_2):x_1,x_2\in{\mathbb Z}\}$ be the set of all 2-vectors of integers.
Elements of $\zt$ are called sites.
The {\dof neighbourhood} $N(x)$ of a site is
$$N(x)=\{y\in\zt: \|x-y\|=1\},$$
where $\|\cdot\|$ is Euclidean distance.  So $|N(x)|=4$ for all $x$.
Let $K$ be a subset of $\zt$.  Define ${\cal B}(K)$ by
$${\cal B}(K)=K\cup\{x\in\zt: |N(x)\cap K|\geq 2\},$$
and $\langle K\rangle$ by
$$\langle K\rangle=\lim_{t\rightarrow\infty} {\cal B}^t(K),$$
where ${\cal B}^t$ denotes the $t$-th iterate.
(If $K$ is the active set at time $0$, then ${\cal B}^t$ is the active set at time $t$ and $\langle K\rangle$ is the active set at time $\infty$). 

Now fix $p\in[0,1]$ and let $X$ be a random subset of $\zt$ in which each site is independently included with probability $p$.  More formally, denote by $P_p$ the product probability measure with parameter $p$ on the product $\sigma$-algebra of $\{0,1\}^{\zt}$, and define the random variable $X$ by $X(\omega)=\{x\in\zt: \omega(x)=1\}$ for $\omega\in\{0,1\}^{\zt}$.  A site $x\in\zt$ is said to be {\dof occupied} if $x\in X$.

We say that a set $K\subseteq\zt$ is {\dof internally spanned} if $\langle X\cap K\rangle=K$.
A {\dof rectangle} is a set of sites of the form
$$R(a,b;c,d):=\{a,\ldots, c\}\times\{b,\ldots, d\},$$
where $a\leq c$, $b\leq d$ are integers.
We also write $R(c,d)=R(1,1;c,d)$.
Define the function
$$I(L,p)=P_p\bigg(R(L,L) \text{ is internally spanned}\bigg).$$
Our main result is the following.
\begin{thm}
\label{main}
Let $L_n,p_n$ be sequences such that $L_n\rightarrow\infty$ and $p_n\rightarrow 0$.
Then:
\begin{mylist}
\item
if $\displaystyle \liminf_{n\rightarrow\infty}\;p_n \log L_n >\lambda$ then $\displaystyle \lim_{n\rightarrow\infty} I(L_n,p_n)=1$;
\item
if $\displaystyle \limsup_{n\rightarrow\infty}\;p_n \log L_n <\lambda$ then $\displaystyle \lim_{n\rightarrow\infty} I(L_n,p_n)=0$;
\end{mylist}
where
$$\lambda=\frac{\pi^2}{18}.$$
\end{thm}
The main step in the proof of Theorem \ref{main} will be Theorem \ref{mainprop} below, concerning the probability of much smaller squares being internally spanned.
\pagebreak

\begin{thm}{ }
\label{mainprop}
{\ }
\begin{mylist}
\item
$\displaystyle \limsup_{p\rightarrow 0} \sup_{m\geq 1} -p\log I(m,p)\leq 2\lambda$.
\item
$\displaystyle \liminf_{B\rightarrow \infty}\liminf_{p\rightarrow 0} -p\log I(\lfloor B/p\rfloor,p) \geq 2\lambda$.
\end{mylist}
\end{thm}
(Here $\lfloor\cdot\rfloor$ denotes the integer part). 

We also establish the following result about the time-evolution of the process on $\zt$.
Define
$$J(t,p)=P_p\bigg((0,0)\in {\cal B}^t(X)\bigg),$$
(the probability the origin is active by time $t$).
\begin{thm}
\label{time}
Let $t_n,p_n$ be sequences such that $t_n\rightarrow\infty$ and $p_n\rightarrow 0$.
Then:
\begin{mylist}
\item
if $\displaystyle \liminf_{n\rightarrow\infty}\;p_n \log t_n >\lambda$ then $\displaystyle \lim_{n\rightarrow\infty} J(t_n,p_n)=1$;
\item
if $\displaystyle \limsup_{n\rightarrow\infty}\;p_n \log t_n <\lambda$ then $\displaystyle \lim_{n\rightarrow\infty} J(t_n,p_n)=0$.
\end{mylist}
\end{thm}

The {\dof modified bootstrap model} is a variant of the above model, in which the definition of ${\cal B}$ is replaced with
$${\cal B}'(K)=K\cup\bigg\{x\in \zt: |\{x+e_i,x-e_i\}\cap K|\geq 1 \text{ for each of }i=1,2\bigg\},$$
where $e_1=(1,0)$ and $e_2=(0,1)$.  We define $I'(\cdot,\cdot)$ and $J'(\cdot,\cdot)$ accordingly.  The arguments in this article can be used to prove the following.
\begin{thm}
\label{modified}
For the modified bootstrap model, the analogues of Theorems \ref{main},\ref{mainprop},\ref{time} hold, with $\lambda$ replaced by
$$\lambda'=\frac{\pi^2}{6}.$$
\end{thm}

Our results answer the main question posed in \cite{aizenman-leb} in the case of two-dimensional bootstrap percolation.  In that paper is was proved (for a wider class of models) that for $p\rightarrow 0$ and $L\rightarrow \infty$, we have $I(L,p)\rightarrow 1$ if $\liminf p\log L>c_1$ and $I(L,p)\rightarrow 0$ if $\limsup p\log L<c_2$, for {\em different} constants $c_1,c_2$, and similarly for $J(L,p)$.  Other aspects of the model have been subsequently studied in detail (see \cite{a-m-s},\cite{gravner-griffeath-2} and the references therein, for example), but the natural question of whether $c_1,c_2$ could be replaced with a single sharp threshold remained open.  Our results answer this affirmatively, as well as establishing the precise value of the threshold.

Predictions for the thresholds $\lambda,\lambda'$ based on simulation are not in good agreement with our rigorous result.  In \cite{adler},\cite{adler2}, the estimates $0.245 \pm 0.015$ for $\lambda$ and $0.47 \pm 0.02$ for $\lambda'$ are reported, whereas $\lambda=\pi^2/18=0.548311\cdots$ and $\lambda'=\pi^2/6=1.644934\cdots$.  The likely reason for this discrepancy is that it is necessary to take $L$ extremely large in order to ``see'' the true limiting behaviour.  (These simulations used values of $L$ up to $28,800$).  Similar phenomena have been noted for several other variants of bootstrap percolation; for details see \cite{frobose},\cite{schonmann-majority},\cite{schonmann},\cite{adler-duarte-enter2},\cite{adler-duarte-enter}.

Bootstrap percolation in three and higher dimensions presents a new set of challenges; for details see \cite{cerf-cirillo},\cite{cerf-manzo},\cite{schonmann}.  In particular for the three-dimensional model in which a site becomes active if at least three of its six neighbours are active, it was established in \cite{cerf-cirillo} that the form of the threshold regime is different: $I(L,p)\rightarrow 1$ if $\liminf p\log \log L>c_1$ and $I(L,p)\rightarrow 0$ if $\limsup p\log \log L<c_2$; the extension to other dimensions is treated in \cite{cerf-manzo}.  The arguments required for these results are more sophisticated than those in \cite{aizenman-leb}, and it seems that the problem of establishing a sharp threshold here is likely to be correspondingly harder.

Aside from their intrinsic mathematical interest, bootstrap percolation models find numerous applications, both directly and as tools in the analysis of more complicated systems.  See for example \cite{adler},\cite{fontes-sidoravicius-schonmann},\cite{frobose},\cite{kirkpatrick},\cite{wolfram}.

In \cite{aizenman-leb} it was also conjectured that stochastic Ising models show similar metastability behaviour.  The analogue of the result in \cite{aizenman-leb} was proved (in arbitrary dimension) in \cite{schonmann-ising}, and the analogue of our result was proved in \cite{schonmann-shlosman}.

We omit the proof of Theorem \ref{modified} regarding the modified bootstrap model.  The proof is almost identical to those given for the earlier theorems, and in fact a few simplifications are possible.  The main difference is that the definition of a rectangle being ``horizontally (respectively vertically) traversable'' in Section \ref{secbasic} should be replaced with the statement that {\em all} the columns (respectively rows) are occupied, and the function $g$ (see Sections \ref{secintegrals}, \ref{secbasic}, \ref{seclower}) should be replaced with $f$.

The article is organized as follows.  In Section \ref{secintegrals} we introduce and solve a definite integral which gives rise to the constant $\lambda=\pi^2/18$.  In Section \ref{secbasic} we give notation and basic results.  In Section \ref{seclower} we prove Theorem \ref{mainprop} (i), and in Section \ref{secmeta} we deduce Theorems \ref{main} and \ref{time} from Theorem \ref{mainprop}.  These proofs are relatively standard; the rest of the article is devoted to the much harder task of proving Theorem \ref{mainprop} (ii).  In Sections \ref{secvar}, \ref{secbord}, \ref{secdisj} and \ref{sechier} we assemble the tools needed for this; the reader may prefer to skip the proofs in these sections on a first reading.  Finally in Section \ref{secupper} we complete the proof of Theorem \ref{mainprop} (ii).

We now present a sketch of the ideas behind the proofs in this article.  Theorem \ref{main} may be deduced from Theorem \ref{mainprop} as follows.  Roughly speaking, Theorem \ref{mainprop} states that if we take $B$ large and then let $p\rightarrow 0$, then $I(B/p,p)$ behaves roughly as
$e^{-2 \lambda_B/p},$
where $\lambda_B\rightarrow \lambda$ as $B\rightarrow \infty$.  Now suppose $L\approx e^{c/p}$ for some $c$.  It may be shown that the event that $R(L,L)$ is internally spanned is roughly equivalent to the event that it contains some internally spanned square of side length $B/p$.  This is because for $B$ large, most intervals of length $B/p$ contain an occupied site, so when a growing active cluster reaches approximately this size it is almost certain to grow to fill the whole square (such a growing cluster is sometimes referred to as a ``critical droplet'' in the literature).  The expected number of internally spanned squares of side $B/p$ in $R(L,L)$ is approximately
$L^2 e^{-2 \lambda_B/p} \approx e^{2(c-\lambda_B)/p},$
which is $\gg 1$ if $c>\lambda_B$ and $\ll 1$ if $c<\lambda_B$.  However, we may choose $B$ so that $\lambda_B$ is as close to $\lambda$ as desired.  A more careful application of these ideas proves Theorem \ref{main}, and a similar approach gives Theorem \ref{time}.

We now turn to the proof of Theorem \ref{mainprop}.  The lower bound on $I$ in (i) is relatively straightforward; the upper bound in (ii) is much harder.  For (i), assume $m=B/p$ and consider one way in which $R(m,m)$ can be internally spanned.  Suppose that $(1,1)$ is occupied, and that for each $k=2,\ldots ,m$, the column with horizontal coordinate $k$ and height $k-1$ contains at least one occupied site, and similarly the row with vertical coordinate $k$ and width $k-1$.  Then it is easily seen that each square with bottom-left corner (1,1) is internally spanned, so the active set grows to fill whole of $R(m,m)$.  This picture gives the correct lower bound for the {\em modified} bootstrap model, but for the bootstrap model one can do better.  Since a site becomes active if it has active neighbours on two opposite sides, we can allow some of the rows and columns described above to be vacant.  Roughly speaking, the growth will still take place provided there are never two adjacent vacant rows or columns (actually, a little more care is needed, see Section \ref{seclower}).  A calculation shows that the probability of this event behaves roughly as
$e^{-(2/p)\int_0^B g}$
where $g$ is a certain function which satisfies $\int_0^\infty g(z) \; dz=\lambda=\pi^2/18$.  This proves Theorem \ref{mainprop} (i), and also shows that in some sense the ``natural'' length unit for the problem is $1/p$.

For the upper bound, Theorem \ref{mainprop} (ii), the basic idea is to consider all other possible ways in which $R=R(B/p,B/p)$ could be internally spanned, and find upper bounds on the number of ways and the probability of each one.  One way in which $R$ could be internally spanned is for the growth described above to start from the center of the square $R$ (say), and for an active square to grow in all four directions rather than just two.  It turns out that the probability of this is essentially the same, $e^{-(2/p)\int_0^B g}$, roughly because growth by one unit each to the left and to the right is equivalent to growth by two units to the right.  Hence there is in some sense no loss of generality in assuming that growth starts from the corner.  However, there are many other ways in which this growth could occur.  For example two adjacent rows as described above might be vacant, in which case the active region could still continue to grow horizontally for a while, until it encounters an occupied site in the vertical direction, at which point vertical growth can resume.  It can be shown that such a possibility has much smaller probability than the one considered earlier.  Indeed, each such growth history corresponds to an oriented path $\gamma$ in $[0,\infty)^2$ (where the length is rescaled by a factor $1/p$), and it may be shown that the probability of such a history is roughly
$e^{-w(\gamma)/p}$
where $w$ is a certain functional on paths defined as a path integral involving the function $g$.  By the convexity of $g$, it can be shown (see Section \ref{secvar}) that this functional is minimized along the main diagonal $x_1=x_2$, with minimum $2\int_0^B g$, which corresponds to growth as a square as described above.

At this point it is natural to try to get an upper bound on the probability of $R$ being internally spanned by summing the above probability over all possible paths.  However, the number of such paths is large, roughly $2^{B/p}$, so this does not give the correct bound.  The solution is to introduce a new ``coarse-graining'' length scale $T/p$, which is small compared with $B/p$ but large compared with the lattice spacing $1$.  It may be shown that the probability of seeing some path which coincides with $\gamma$ at a set of points of spacing $T$ is also approximately
$e^{-w(\gamma)/p}.$
This is because all such paths require roughly the same occupied sites.  Actually, for this argument to work, we must assume that the growth starts not from the point $(1,1)$ but from a larger rectangle, of size $A/p$, where $T\ll A\ll 1 \ll B$; the details are in Section \ref{secbord}.  The number of possible choices of a set of points at spacing $T/p$ is much less than before (in fact, bounded in $p$), so now summing over all possibilities gives an upper bound of approximately
$e^{-(2/p)\int_A^B g}$
(and the integral here is close to $\int_0^\infty g$ for $A$ small and $B$ large, since the latter integral converges).

However, there is another difficulty.  There are other ways for $R$ to be internally spanned which do not involve only growth starting from a single ``seed''.  For example, the growing square described above might encounter another internally spanned rectangle $S$ within $R$, and the two would combine to give a larger internally spanned rectangle without any of the intervening growth taking place.  It seems unlikely that such an event could have higher probability, since the probability of a small internally spanned rectangle $S$ should be much less than the probability of finding at least one occupied site in each of the corresponding rows and columns required for growth as described previously.  In fact, most of the work in the proof is to rule out possibilities such as this.

The idea is as follows.  Suppose $R$ is internally spanned.  Then we expect to be able to find a sequence of successively slightly smaller internally spanned rectangles inside $R$ (corresponding to the growth picture above).  However, at some point there may be a ``split'' into two separate internally spanned rectangles.  Now we apply the same reasoning to each of these, and so on.  In this way we obtain a ``hierarchy'' - a tree structure of nested rectangles each of which is internally spanned.  We can arrange that when there is no splitting, the sizes of consecutive rectangles differ by approximately the coarse-graining scale $T/p$; when there is a split, the two offspring rectangles might be much smaller than the parent, but they should have the property that they ``span'' the parent, in the sense that bootstrap percolation starting with the two offspring rectangles completely filled results in the parent being completely filled.  Also, for the coarse-graining argument to work, we should stop whenever the rectangles in a line of descent become too small (perhaps smaller than $A/p$), and declare such a rectangle to be a seed.

Now, another difficulty arises.  Since there may be many splits, there may also be many seeds; indeed $R$ could be ``almost all seeds'' in which case there is no room left for the growth estimates used previously.  The solution is to use the fact that seeds are small, together with an a priori bound on the probability a small rectangle is internally spanned, to show that there are ``not too many'' seeds.  More precisely, we introduce yet another length $Z/p\ll A/p$, and declare a rectangle a seed if its size is less than $Z/p$ (the exact meaning of ``size'' turns out to be a delicate issue here).  We show that the probability that the total perimeter of all the seeds is more than $A/p$ is very small.  Then, provided the total perimeter is less than $A/p$, a geometrical argument together with a variational principle can be used to show that the probability of the whole hierarchy is approximately at most
$e^{-(2/p)\int_A^B g}.$
The order of choosing the constants is important here: they must be chosen in the order $B,A,Z,T$, then finally $p\rightarrow 0$.  The total number of hierarchies is large, but it can be shown to be at most
$p^{-K}$
where $K$ depends on $B$ but not $p$.  Hence summing gives a bound of the correct form.

At the heart of any such proof must be a (deterministic) result giving rigorous necessary conditions for a rectangle to be internally spanned.  In \cite{aizenman-leb}, this role was played by a result which we shall also make use of, Lemma \ref{al}.
For the construction of hierarchies, the required condition is provided by an apparently new result, Proposition \ref{disjoint}, which roughly speaking states that any internally spanned rectangle must be spanned by a pair of smaller rectangles which are themselves internally spanned {\em disjointly} (in the sense of the Van den Berg-Kesten inequality).

\section{Integrals}
\label{secintegrals}

We define the functions $f:(0,\infty)\rightarrow (0,\infty)$, $\beta:(0,1)\rightarrow (0,1)$, and  $g:(0,\infty)\rightarrow (0,\infty)$ by
$$f(z)=-\log(1-e^{-z}),$$
$$\beta(u)=\frac{u+\sqrt{u(4-3u)}}{2},$$
$$g(z)=-\log\beta(1-e^{-z}).$$
Note that $f$ and $g$ are continuously differentiable, positive, decreasing and convex.  To see the latter for $g$ note that the function $-z \mapsto g(z)$ is the composition of the increasing, convex functions $z\mapsto e^z$, $z\mapsto -\beta(1-z)$ and $z\mapsto -\log(-z)$.  Observe also that $\beta(u)>u$ so $f(z)>g(z)$ for $z\in (0,\infty)$.  Finally note that $f(z),g(z) \rightarrow \infty$ as $z\rightarrow 0$, and $f(z),g(z),zf(z),zg(z)\rightarrow 0$ as $z\rightarrow \infty$.

\begin{prop} \nopagebreak
\label{integrals}
{\ }
\begin{mylist}
\item
$$\int_0^\infty f(z)\; dz=\frac{\pi^2}{6};$$ \nopagebreak
\item
$$\int_0^\infty g(z)\; dz=\frac{\pi^2}{18}.$$
\end{mylist}
\end{prop}
We define $\lambda$ to equal the second integral.

Proposition \ref{integrals} (ii) is somewhat remarkable; the integral does not appear in standard tables such as \cite{grad}, and seems not to be directly solvable by computer programs such as Mathematica.  The proof which we shall give is mysterious, and uses many special features of the function $g$.
\begin{pfof}{Proposition \ref{integrals}}
In (i) we substitute $x=e^{-z}$ to obtain
$$\int_0^\infty -\log(1-e^{-z}) \; dz=\int_0^1 -\log(1-x)\frac{dx}{x}=\frac{\pi^2}{6},$$
by \cite{grad}, number 4.291.2.

The integral in (ii) converges because $g\leq f$.  Substituting $y=e^{-z}$ and then taking out a factor $1-y$, we have
\begin{eqnarray*}
\lefteqn{\int_0^\infty g(z)\; dz} \\
&=&\int_0^\infty -\log\frac{1-e^{-z}+\sqrt{(1-e^{-z})(1+3e^{-z})}}{2}\;dz \\
&=&\int_0^1 -\log\frac{1-y+\sqrt{(1-y)(1+3y)}}{2}\; \frac{dy}{y} \\
&=&\int_0^1 -\log(1-y)\;\frac{dy}{y}+\int_0^1 
-\log\frac{1+\sqrt{\frac{1+3y}{1-y}}}{2}\; \frac{dy}{y}.
\end{eqnarray*}
Note that the first of these two integrals converges by the above, hence the second must also.

In the first integral we make the substitution $x=1-y$.  In the second we make the substitution
$$x=\frac{2}{1+\sqrt{\frac{1+3y}{1-y}}}$$
which yields
$$y=\frac{1-x}{1-x+x^2}$$
and
$$\frac{dy}{y}=\frac{(2-x)x}{(x-1)(1-x+x^2)} dx,$$
and also interchanges the limits $0$ and $1$.
Thus we obtain
\begin{eqnarray*}
\lefteqn{\int_0^\infty g(z)\; dz} \\
&=&\int_0^1 -\log x\;\frac{dx}{1-x}+\int_0^1  
-\log x \;\frac{(2-x)x}{(x-1)(1-x+x^2)}\; dx \\
&=&\int_0^1 \log x\left(-\frac{1}{1-x}-\frac{(2-x)x}{(x-1)(1-x+x^2)}\right) \;dx \\
&=&\int_0^1 \log x\; \frac{2x-1}{1-x+x^2} \;dx \\
&=&\bigg[ \log x \;\log (1-x+x^2) \bigg]_0^1 -\int_0^1\log(1-x+x^2) \;\frac{dx}{x} \quad\text{(by parts)} \\
&=&0-\int_0^1 \log\left(\frac{1+x^3}{1+x}\right) \;\frac{dx}{x} \\
&=&\int_0^1 \log(1+x)\;\frac{dx}{x} -\int_0^1 \log(1+x^3)\;\frac{dx}{x} \\
&=&\int_0^1 \log(1+x)\;\frac{dx}{x} -\int_0^1 \log(1+u)\;\frac{du}{3u} \quad\text{(substituting $u=x^3$)}\\
&=&\frac{2}{3}\int_0^1 \log(1+x)\;\frac{dx}{x}\\
&=& \frac{2}{3} \;\frac{\pi^2}{12} \qquad\text{(by \cite{grad}, 4.291.1)} \\
&=& \frac{\pi^2}{18}.
\end{eqnarray*}
\end{pfof}

\section{Basic Notation and Results}
\label{secbasic}

In this section we introduce basic tools which will be used throughout.
The following elementary properties of $\langle\cdot\rangle$ will be fundamental: $\langle K\rangle \supseteq K$; $\langle\langle K\rangle\rangle = \langle K\rangle$; if $K \subseteq K'$ then $\langle K\rangle \subseteq \langle K'\rangle$.  These have the important consequence that if $K\subseteq K'\subseteq \langle K\rangle$ then $\langle K'\rangle = \langle K\rangle$.

Let $R=R(a,b;c,d)$ be a rectangle.
By the {\dof dimensions} of $R$ we mean the 2-vector
$$\dim(R)=(c-a+1,d-b+1).$$
If $\dim(R)=(m,n)$ we define the {\dof short side} $\sht(R)=\min\{m,n\}$, the {\dof long side} $\lng(R)=\max\{m,n\}$ and the {\dof semi-perimeter} $\phi(R)=m+n$ of $R$.

A site $x\in\zt$ is {\dof occupied} if $x\in X$.  A set of sites $K\subseteq\zt$ is {\dof full} if every site in $K$ is occupied, and {\dof occupied} if at least one site in $K$ is occupied.
It will be convenient to write
$$q=-\log(1-p).$$
Note that $q\geq p$ and $q \sim p$ (that is, $q/p \rightarrow 1$) as $p\rightarrow 0$.  The advantage of this notation is that
\begin{equation}
\label{occ}
P_p(K \text{ is occupied})=1-e^{-|K|q}=e^{-f(|K|q)}
\end{equation}
(where $f$ was defined in Section \ref{secintegrals}).

We say that a sequence of events $(A_j)$ has a {\dof double gap} if there is an adjacent pair $A_i,A_{i+1}$ neither of which occurs.

\begin{lemma}
\label{doublegap}
In a sequence of $k$ independent events each with probability $u\in(0,1)$, the probability $a_k(u)$ that there are no double gaps satisfies
$$\beta(u)^k\leq a_k(u)\leq \beta(u)^{k-1}$$
where
$$\beta(u)=\frac{u+\sqrt{u(4-3u)}}{2}.$$
\end{lemma}

\begin{pf}
By induction on $k$, on noting that $a_0=a_1=1$, $a_{k+2}=ua_{k+1}+(1-u)ua_k$, $0<\beta<1$, and $\beta^2=u\beta+(1-u)u$.
\end{pf}

Let $R=R(a,b;c,d)$.  For $a\leq i\leq c$ (respectively $b\leq j\leq d$), {\dof column} $i$ (respectively {\dof row} $j$) of $R$ is the rectangle $R(i,b;i,d)$ (respectively $R(a,j;c,j)$).  We say that $R$ is {\dof horizontally} (respectively {\dof vertically}) {\dof traversable} if the sequence $(\text{column $i$ is occupied})_{i=a}^c$ (respectively $(\text{row $j$ is occupied})_{j=b}^d$) has no double gaps.  The following definitions will also be convenient.  The rectangle $R$ is {\dof East-traversable} (respectively {\dof North-traversable}) if it is horizontally (respectively vertically) traversable and in addition column $c$ (respectively row $d$) is occupied.

Recall the definitions of $f$ and $g$ from Section \ref{secintegrals}.
\begin{lemma}
\label{traversable}
If $R$ is a rectangle with dimensions $(m,n)$ then
\begin{mylist}
\item
$$e^{-mg(nq)} \leq P_p(R \text{ is horizontally traversable})\leq e^{-(m-1)g(nq)};$$
\item
$$e^{-mf(nq)}\leq e^{-(m-1)g(nq)-f(nq)} \leq P_p(R \text{ is East-traversable})\leq e^{-mg(nq)}.$$
\end{mylist}
And similar inequalities hold for vertical and North- traversability, with $m$ and $n$ exchanged.
\end{lemma}

\begin{pf}
Part (i) follows immediately from (\ref{occ}), Lemma \ref{doublegap} and the definition of $g$.  For the upper bound in (ii), note that if $R(a,b;c,d)$ is East-traversable then $R(a,b;c+1,d)$ is horizontally traversable, and use (i).  The second inequality in (ii) is straightforward, and the first follows because $f\geq g$.
\end{pf}

To see the usefulness of the above concepts note the following.
\begin{lemma}
\label{travapps}
{\ }
\begin{mylist}
\item
If $R$ is internally spanned then $R$ is East- and North- traversable.
\item
If $R_1=(a,b;c,d)$ is internally spanned and $R_2=(c+1,b;e,d)$ is East-traversable then $R_1\cup R_2$ is internally spanned.  And a similar statement holds for North-traversability.
\end{mylist}
\end{lemma}

\begin{pf}
For part (i), note that if two adjacent columns or the East-most column of $R$ contains no occupied sites, then no site in these columns can be in $\langle R\cap X\rangle$.  A similar remark applies to rows.  For (ii), it is easy to see that each successive column in $R_2$ (moving in the East direction) is in $\langle R\cap X\rangle$.
\end{pf}

An event $A$ of $\{0,1\}^{\zt}$ is called {\dof increasing} if whenever $\omega\in A$ and $\omega'\geq \omega$ we have $\omega'\in A$.  The FKG inequality states that for increasing events $A,B$ we have $P_p(A\cap B)\geq P_p(A)P_p(B)$ (see \cite{g2} p34 for example).

\section{Lower Bound}
\label{seclower}

\begin{pfof}{Theorem \ref{mainprop} (i)}
Let $r=\lfloor p^{-1/2}\rfloor$, and let $A=A(m,p)$ be the event that all the following occur:
$$R(1,1;1,r) \text{ and } R(1,1;r,1) \text{ are full},$$
$$\text{the sites } (m,1) \text{ and } (1,m) \text{ are occupied},$$
\hspace{\parindent} and for all integers $k\geq 1$:
$$R(kr+1,1;kr+r,kr) \text{ is East-traversable},$$
$$R(1,kr+1;kr,kr+r) \text{ is North-traversable}.$$
See Figure \ref{defa}.
\begin{figure}
\centering
\begin{picture}(0,0)%
\includegraphics{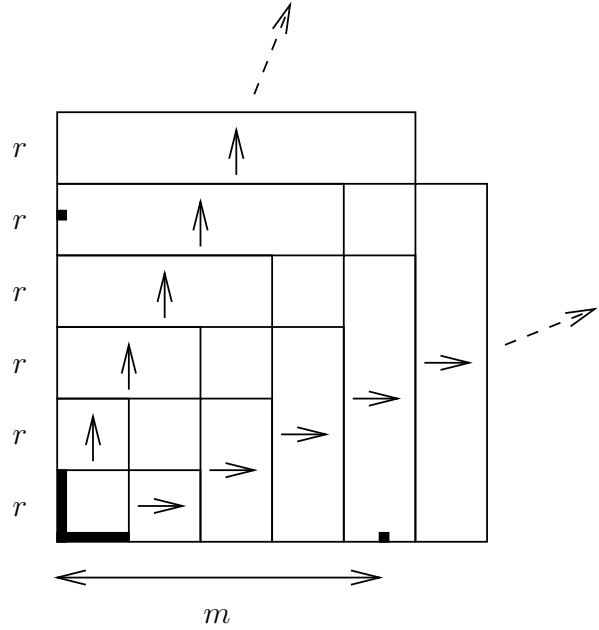}%
\end{picture}%
\setlength{\unitlength}{0.00062500in}%
\begin{picture}(4876,5248)(547,-4687)
\put(601,-3736){\makebox(0,0)[b]{$r$}}
\put(601,-3136){\makebox(0,0)[b]{$r$}}
\put(601,-2536){\makebox(0,0)[b]{$r$}}
\put(601,-1936){\makebox(0,0)[b]{$r$}}
\put(601,-1336){\makebox(0,0)[b]{$r$}}
\put(601,-736){\makebox(0,0)[b]{$r$}}
\put(2251,-4636){\makebox(0,0)[b]{$m$}}
\end{picture}
\caption{An illustration of the event $A$.  The arrows indicate East- and North- traversability.}
\label{defa}
\end{figure}

Using Lemma \ref{travapps} (ii), it is easily seen that if $A$ occurs then $R(m,m)$ is internally spanned, hence using Lemma \ref{traversable} (ii) and the FKG inequality,
$$I(m,p)\geq P_p(A)\geq
p^{2r+1}\prod_{k=1}^\infty \bigg(e^{-(r-1)g(krq)-f(krq)}\bigg)^2,$$
for all $m$.
Hence,
\begin{eqnarray*}
\lefteqn{\sup_{m\geq 1} -p\log I(m,p)} \\
 &\leq&
 -(2r+1)p\log p + \frac{2(r-1)p}{rq}\sum_{k=1}^\infty g(krq)rq
 + \frac{2p}{rq}\sum_{k=1}^\infty f(krq)rq\\
&\leq& -3rp\log p + \frac{2p}{q}\int_0^\infty g(z)\;dz + \frac{2p}{rq} \int_0^\infty f(z)\;dz
\end{eqnarray*}
where in the last step we have use the fact that $f$ and $g$ are decreasing.
Now, recalling that $q \sim p$, $r=\lfloor p^{-1/2}\rfloor$, and $\int_0^\infty f(z)\;dz<\infty$, we see that as $p\rightarrow 0$ the above expression converges to $0+2\int_0^\infty g(z)\;dz+0=2\lambda$, as required.
\end{pfof}

\section{Metastability}
\label{secmeta}

In this section we deduce the (i) parts of Theorems \ref{main} and \ref{time} from the (i) part of Theorem \ref{mainprop}, and similarly for the (ii) parts.

\begin{pfof}{Theorem \ref{main} (i)}
It is clearly sufficient to prove that for any $\epsilon>0$, if $p_n\rightarrow 0$ and $L_n\rightarrow\infty$ are such that $p_n\log L_n \geq \lambda+\epsilon$ then $I(L_n,p_n)\rightarrow 1$.

Suppose $p\log L\geq \lambda+\epsilon$ and $p<1/2$.  Let $m=\lfloor p^{-3}\rfloor$, and let $S=S(L,p)$ be the event that $R(L,L)$ contains at least one internally spanned rectangle of dimensions $(m,m)$.  By dividing $R(L,L)$ into disjoint rectangles of dimensions $(m,m)$ we see that
$$P_p(S)\geq 1-\bigg( 1-I(m,p)\bigg)^{\lfloor L/m\rfloor ^2},$$
so
$$-\log(1-P_p(S))\geq \frac{L^2p^6}{2} I(m,p).$$
Therefore
$$p\log\bigg[-\log(1-P_p(S))\bigg]\geq 2p\log L +6p\log p -p\log 2 +p\log I(m,p).$$
Hence, by Theorem \ref{mainprop} (i) we have
$$\liminf_{n\rightarrow\infty} p_n\log\bigg[-\log(1-P_{p_n}(S_n))\bigg]\geq 2(\lambda+\epsilon)+0-0-2\lambda=2\epsilon,$$
where $S_n=S(L_n,p_n)$, so in particular $P_{p_n}(S_n)\rightarrow 1$.

The following is proved in \cite{aizenman-leb}.  There exists a sequence of increasing events $H_n$ defined in terms of the states of sites in $R(L_n,L_n)$ such that $P_{p_n}(H_n)\rightarrow 1$, and if $S_n$ and $H_n$ occur then $R(L_n,L_n)$ is internally spanned.  This proves the result, since by the FKG inequality,
$$I(L_n,p_n)\geq P_{p_n}(S_n)P_{p_n}(H_n)\rightarrow 1.$$

Here is a sketch of the construction of $H_n$.  If $p_n\log L_n\rightarrow a>\lambda$ then we may take $H_n$ to be the event that every rectangle of dimensions $(m,1)$ or $(1,m)$ in $R(L_n,L_n)$ is occupied.  If $L_n$ grows faster than this then we must also use a renormalization argument, dividing $R(L_n,L_n)$ into disjoint squares of size $e^{(\lambda+\delta)/p_n}$, so that the probability each one is internally spanned exceeds the critical probability for site percolation on $\zt$.
\end{pfof}

\begin{pfof}{Theorem \ref{time} (i)}
Note that, in contrast with $I$, $J(t,p)$ is clearly increasing in $t$, so we may assume that $p_n\log t_n\rightarrow \lambda +\epsilon$, where $\epsilon>0$.  Let $m_n=\lfloor p_n^{-3}\rfloor$ and $L_n=\lfloor t_n/(3m_n)\rfloor$; then $p_n\log L_n\rightarrow \lambda +\epsilon$.  It is easily seen that the events $H_n$ in the proof of Theorem \ref{main} (i) above may be chosen in such a way that if $S_n$ and $H_n$ occur then the whole of $R(L_n,L_n)$ becomes active in time at most
$$|R(m_n,m_n)| + 2m_nL_n \leq t_n$$
(The inequality holds if $p_n$ is sufficiently small).  Hence,
$$P_n\bigg( (1,1) \in {\cal B}^{t_n} (X\cap R(L_n,L_n)) \bigg)\rightarrow 1,$$
therefore $J(t_n,p_n)\rightarrow 1$.
\end{pfof}

To prove Theorem \ref{main} (ii) we need the following result from \cite{aizenman-leb}.
\begin{lemma}
\label{al}
Let $k$ be a positive integer and let $R$ be a rectangle with $\lng(R)$ $\geq 2k+1$.  If $R$ is internally spanned then there exists an internally spanned rectangle $T\subseteq R$ with $\lng(T)\in[k,2k+1].$
\end{lemma}
The idea of the proof of Lemma \ref{al} is to construct the bootstrap percolation process by an algorithm which sequentially replaces a pair of internally spanned rectangles with a larger internally spanned rectangle; the result then follows because at each step the long side of the longest rectangle increases at most as $n\mapsto 2n+1$.  For the details see \cite{aizenman-leb}.  We shall use a similar argument in Section \ref{secdisj} to prove Proposition \ref{disjoint}.

\begin{pfof}{Theorem \ref{main} (ii)}
It is clearly sufficient to show that for any $\epsilon>0$, if $p_n\rightarrow 0$ and $L_n\rightarrow\infty$ with $p_n\log L_n\leq\lambda-\epsilon$ then $I(L_n,p_n)\rightarrow 0$.
We write $\lambda_B=\frac{1}{2}\inf_{B'\geq B}\liminf_{p\rightarrow 0} -p\log I(\lfloor B'/p\rfloor,p)$, so that Theorem \ref{mainprop} (ii) states $\lim_{B\rightarrow \infty} \lambda_B\geq \lambda$.

Fix $\epsilon>0$ and $B>0$, and suppose that $p\log L\leq \lambda-\epsilon$. We write $R=R(L,L)$ and let $K=\lfloor B/p\rfloor$ and $k=\lfloor B/(2p)\rfloor-1$, so that $k\leq 2k+1\leq K$.  Assume that $p$ is sufficiently small that $1<k<K<L$.
 We claim first that if $T$ is any rectangle satisfying $\lng(T)\in [k,2k+1]$ then
\begin{equation}
\label{compare}
P_p(T \text{ is internally spanned})\leq e^{2Kf(kq)} I(K,p).
\end{equation}
To prove this, suppose without loss of generality that $T$ is of the form $T=R(a,b)$ where $a\in [k,2k+1]$ and $b\leq K$.  Let $S_1=R(1,b+1;a,K)$ and $S_2=R(a+1,1;K,K)$; see Figure \ref{s1s2}.  By Lemma \ref{travapps} (ii), if $T$ is internally spanned, and $S_1$ is North-traversable, and $S_2$ is East-traversable, then $R(K,K)$ is internally spanned.  Hence by Lemma \ref{traversable} (ii) we have
$$I(K,p)\geq P_p(T \text{ is internally spanned}) e^{-(K-b)f(aq)-(K-a)f(Kq)}$$
which yields (\ref{compare}) since $f$ is decreasing.
\begin{figure}
\centering
\begin{picture}(0,0)%
\includegraphics{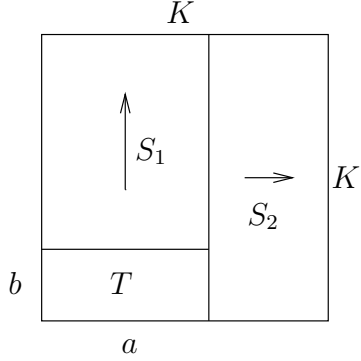}%
\end{picture}%
\setlength{\unitlength}{0.00041700in}%
\begin{picture}(4370,4383)(831,-4363)
\put(2251,-4336){{$a$}}
\put(901,-3586){\makebox(0,0)[b]{$b$}}
\put(3001,-136){\makebox(0,0)[b]{$K$}}
\put(5101,-2236){\makebox(0,0)[b]{$K$}}
\put(2251,-3586){\makebox(0,0)[b]{$T$}}
\put(2626,-1936){\makebox(0,0)[b]{$S_1$}}
\put(4051,-2761){\makebox(0,0)[b]{$S_2$}}
\end{picture}
\caption{The rectangles $T,S_1,S_2$.  The arrows indicate North- and East- traversability.}
\label{s1s2}
\end{figure}

Hence using Lemma \ref{al} and (\ref{compare}) we have
$$I(L,p)\leq L^2K^2 e^{2Kf(kq)} I(K,p).$$
Here $L^2K^2$ is an upper bound on the number of possible choices for the rectangle $T$.
Hence
$$p\log I(L,p) \leq 2p\log L +2p\log K+2Kpf(kq) +p\log I(K,p).$$
Therefore, recalling that $p_n\log L_n\leq\lambda-\epsilon$, and the definitions of $k,K,\lambda_B$, we have 
\begin{eqnarray*}
\limsup_{n\rightarrow \infty} p_n\log I(L_n,p_n) &\leq &
 2(\lambda-\epsilon)+0+2Bf(B/2)-2\lambda_B.
\end{eqnarray*}
Since $\lim_{B\rightarrow\infty} Bf(B/2)= 0$ and $\lim_{B\rightarrow \infty} \lambda_B\geq \lambda$, we may choose $B$ sufficiently large (depending on $\epsilon$) that the above expression is negative, so in particular $I(L_n,p_n)\rightarrow 0$.
\end{pfof}

\begin{pfof}{Theorem \ref{time} (ii)}
In consequence of the results (5.4),(5.5) of \cite{aizenman-leb}, if $I(L,p)\rightarrow 0$ then $J(L,p)\rightarrow 0$.  Hence the result follows from Theorem \ref{main} (ii).
\end{pfof}

\section{Variational Principles}
\label{secvar}

We write $\rtp=(0,\infty)^2=\{\vc{a}=(\vc{a}_1,\vc{a}_2):\vc{a}_1,\vc{a}_2\in (0,\infty)\}$ for the set of all 2-vectors of positive reals.  Boldface letters will denote elements of $\rtp$ unless stated otherwise.
We write $\vc{a}\leq \vc{b}$ if $\vc{a}_1\leq \vc{b}_1$ and $\vc{a}_2\leq \vc{b}_2$.  Let $g$ be a continuously differentiable, positive, decreasing, convex function from $(0,\infty)$ to $(0,\infty)$ ($g$ is otherwise arbitrary for the purposes of this section).  A (piecewise-linear, oriented) {\dof path} $\gamma$ from $\vc{a}$ to $\vc{b}$ is a subset of $\rtp$ consisting of the union of a finite sequence of points $\vc{a}=\vc{u}^0\leq \vc{u}^1\leq\cdots\leq \vc{u}^k=\vc{b}$ (called the {\dof vertices} of the path) together with the line segments $\{\alpha\vc{u}^i+(1-\alpha)\vc{u}^{i+1}:\alpha\in(0,1)\}$ joining the vertices in order.  A path $\gamma$ may be parameterized as $\gamma=\{(x(t),y(t)):t\in [a,b]\}$ where $x(t)+y(t)=t$.  We define the functional $w(\gamma)$ as the path integral
$$w(\gamma)=\int_\gamma \bigg( g(y)dx+g(x)dy \bigg)
=\int_a^b \left(g(y)\frac{dx}{dt}+g(x)\frac{dy}{dt}\right)dt,$$
and for $\vc{a}\leq \vc{b}$ we define
$$W(\vc{a},\vc{b})=\inf_{\gamma:\vc{a}\rightarrow\vc{b}} w(\gamma)$$
where $\inf_{\gamma:\vc{a}\rightarrow\vc{b}}$ denotes the infimum over all paths $\gamma$ from $\vc{a}$ to $\vc{b}$.

The purpose of this self-contained section is to prove the following four properties of $W$.  Propositions \ref{vsum}, \ref{vupper} and \ref{vexact} are fairly natural; Proposition \ref{vsplit} is tailored to a specific application in Section \ref{secupper}.
\begin{prop}
\label{vsum}
If $\vc{a}\leq\vc{b}\leq\vc{c}$ then
$W(\vc{a},\vc{b})+W(\vc{b},\vc{c})\geq W(\vc{a},\vc{c})$.
\end{prop}
\begin{prop}
\label{vupper}
$W(\vc{a},\vc{b})\leq (\vc{b}_1-\vc{a}_1)g(\vc{a}_2)+(\vc{b}_2-\vc{a}_2)g(\vc{a}_1)$.
\end{prop}
\begin{prop}
\label{vexact}
If $\vc{a}_1+\vc{a}_2=A$ and $\vc{b}=(B,B)$ where $A\leq B$ then
$$W(\vc{a},\vc{b})\geq 2 \int_A^B g(z)\; dz.$$
\end{prop}
\begin{prop}
\label{vsplit}
Suppose that $\vc{a}\leq\vc{b}$; $\vc{c}\leq\vc{d}$; $\vc{r}\geq\vc{b}$; $\vc{r}\geq\vc{d}$; $\vc{r}\leq\vc{b}+\vc{d}+(q,q)$ and $\vc{r}\geq (2Z,2Z)$ where $q<Z$.  Then there exists $\vc{s}$ satisfying $\vc{s}\leq \vc{r}$ and $\vc{s}\leq \vc{a}+\vc{c}$ such that
$$W(\vc{a},\vc{b})+W(\vc{c},\vc{d})\geq W(\vc{s},\vc{r})-2qg(Z).$$
Indeed we may take $\vc{s}=\vc{a}\vee[(\vc{a}+\vc{c})\wedge (\vc{a}+\vc{c}+\vc{r}-\vc{b}-\vc{d})]$.
\end{prop}
In the above, $\vee$ and $\wedge$ denote coordinate-wise maximum and minimum respectively.

\begin{pfof}{Proposition \ref{vsum}}
We have
\begin{eqnarray*}
W(\vc{a},\vc{c})
&=&\inf_{\gamma:\vc{a}\rightarrow\vc{c}} w(\gamma) \\
&\leq&\inf_{\gamma:\vc{a}\rightarrow\vc{b}\rightarrow\vc{c}} w(\gamma) \\
&=&\inf_{\gamma_1:\vc{a}\rightarrow\vc{b}} w(\gamma_1)
  +\inf_{\gamma_2:\vc{b}\rightarrow\vc{c}} w(\gamma_2) \\
&=&W(\vc{a},\vc{b})+W(\vc{b},\vc{c}),
\end{eqnarray*}
where $\inf_{\gamma:\vc{a}\rightarrow\vc{b}\rightarrow\vc{c}}$ denotes the infimum of all paths from $\vc{a}$ to $\vc{c}$ containing $\vc{b}$.
\end{pfof}

\begin{pfof}{Proposition \ref{vupper}}
For any path $\gamma$ from $\vc{a}$ to $\vc{c}$, since $g$ is decreasing and $\vc{u}\geq\vc{a}$ for all $\vc{u}\in\gamma$, we have
\begin{eqnarray*}
w(\gamma)&=&\int_\gamma \bigg( g(y)dx+g(x)dy \bigg)\\
&\leq&\int_\gamma \bigg( g(\vc{a}_2)dx+g(\vc{a}_1)dy \bigg)\\
&=&(\vc{b}_1-\vc{a}_1)g(\vc{a}_2)+(\vc{b}_2-\vc{a}_2)g(\vc{a}_1).
\end{eqnarray*}
\end{pfof}

In order to prove Proposition \ref{vexact} we need Lemma \ref{vdiag} below.
For sets $A,B\subseteq\rtp$ we say $A$ {\dof lies Northwest} of $B$ and write $A\succeq B$ if for any $\vc{a}\in A$ and $\vc{b}\in B$ satisfying $\vc{a}_1+\vc{a}_2=\vc{b}_1+\vc{b}_2$, we have $\vc{a}_2\geq \vc{b}_2$.  Let $\Delta$ be the ``main diagonal'' of $\rtp$:
$$\Delta:=\{\vc{u}\in\rtp: \vc{u}_1=\vc{u}_2\}.$$

\begin{lemma}
\label{vdiag}
If $\gamma_1,\gamma_2$ are paths from $\vc{a}$ to $\vc{b}$, and either $\gamma_1\succeq\gamma_2\succeq\Delta$ or $\Delta\succeq\gamma_2\succeq\gamma_1$, then
$w(\gamma_1)\geq w(\gamma_2)$.
\end{lemma}

\begin{pf}
Without loss of generality (since the definition of $w$ is symmetric in the two coordinates), we may assume $\gamma_1\succeq\gamma_2\succeq\Delta$.

\begin{figure}
\centering
\begin{picture}(0,0)%
\includegraphics{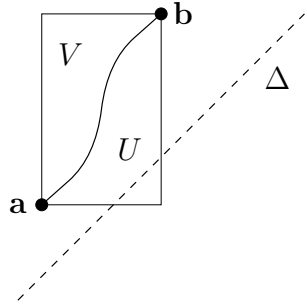}%
\end{picture}%
\setlength{\unitlength}{0.00041700in}%
\begin{picture}(3684,3705)(1139,-3983)
\put(3301,-436){\makebox(0,0)[b]{$\vc{b}$}}
\put(1201,-2836){\makebox(0,0)[b]{$\vc{a}$}}
\put(4501,-1261){\makebox(0,0)[b]{$\Delta$}}
\put(1876,-961){\makebox(0,0)[b]{$V$}}
\put(2626,-2161){\makebox(0,0)[b]{$U$}}
\end{picture}
\caption{The regions $U$ and $V$.}
\label{uv}
\end{figure}
First note that for any path $\gamma$ from $\vc{a}$ to $\vc{b}$, we may express $w(\gamma)$ as
\begin{equation}
\label{areaint}
w(\gamma)=\intint_U g'(y) \;dx\; dy + (\vc{b}_1-\vc{a}_1)g(\vc{a}_2)
+ \intint_V g'(x) \;dx\; dy + (\vc{b}_2-\vc{a}_2)g(\vc{a}_1)
\end{equation}
where the regions $U,V$ are defined by
\begin{eqnarray*}
U&=&\{\vc{u} : \vc{a}\leq\vc{u}\leq\vc{b} \text{ and } \gamma\succeq\{\vc{u}\}\}, \\
V&=&\{\vc{u} : \vc{a}\leq\vc{u}\leq\vc{b} \text{ and } \{\vc{u}\}\succeq\gamma\};
\end{eqnarray*}
see Figure \ref{uv}.  To check (\ref{areaint}) note that by performing the $y$ integral we have
$$ \intint_U g'(y) \;dx\; dy = \int_\gamma \bigg(g(y)-g(\vc{a}_2)\bigg)\; dx = \int_\gamma g(y)\ dx - (\vc{b}_1-\vc{a}_1)g(\vc{a}_2),$$
and similarly for the second integral.

Applying (\ref{areaint}) to $\gamma_1$ and $\gamma_2$ and subtracting, we have 
$$w(\gamma_1)-w(\gamma_2)=\intint_H \bigg( g'(y)-g'(x)\bigg)\;dx\;dy,$$
where $H$ is the region between $\gamma_1$ and $\gamma_2$:
$$H=\{\vc{u} : \vc{a}\leq\vc{u}\leq\vc{b} \text{ and } \gamma_1\succeq\{\vc{u}\}\succeq\gamma_2\}.$$
Since $\gamma_1\succeq\gamma_2\succeq\Delta$, we have $H\succeq\Delta$, and hence $y\geq x$ for $(x,y)\in H$, and since $g$ is convex this implies that $g'(y)-g'(x)\geq 0$ on $H$, so $w(\gamma_1)-w(\gamma_2) \geq 0$.
\end{pf}

\begin{pfof}{Proposition \ref{vexact}}
Let $\gamma$ be a path from $\vc{a}$ to $\vc{b}$.  We claim that
$$w(\gamma)\geq w(\gamma_0)$$ 
where $\gamma_0$ is the path with vertices $\vc{a},\vc{u},\vc{b}$, where $\vc{u}_1=\vc{u}_2=\max\{\vc{a}_1,\vc{a}_2\}$ (see Figure \ref{opt}).  Thus $W(\vc{a},\vc{b})=w(\gamma_0)$.  To check the above claim, split $\gamma$ into sections separated by the intersections of $\gamma$ and $\gamma_0$, and observe that by Lemma \ref{vdiag}, each section of $\gamma$ has a value of $w$ at least as large as the corresponding section of $\gamma_0$.

\begin{figure}
\centering
\begin{picture}(0,0)%
\includegraphics{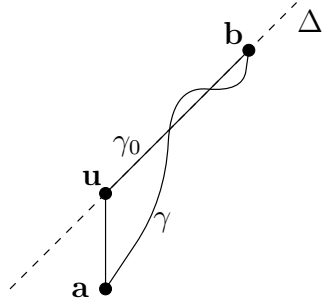}%
\end{picture}%
\setlength{\unitlength}{0.00041700in}%
\begin{picture}(3947,3724)(1479,-3763)
\put(2401,-3736){\makebox(0,0)[b]{$\vc{a}$}}
\put(2551,-2311){\makebox(0,0)[b]{$\vc{u}$}}
\put(5326,-361){\makebox(0,0)[b]{$\Delta$}}
\put(4351,-511){\makebox(0,0)[b]{$\vc{b}$}}
\put(3451,-2911){\makebox(0,0)[b]{$\gamma$}}
\put(3001,-1936){\makebox(0,0)[b]{$\gamma_0$}}
\end{picture}
\caption{A path $\gamma$ from $\vc{a}$ to $\vc{b}$, and the optimizing path $\gamma_0$.}
\label{opt}
\end{figure}

Now since $g$ is positive, $w(\gamma_0)\geq w(\gamma_0')$ where $\gamma_0'$ is the straight path with vertices $\vc{u}$ and $\vc{b}$.  From the definition of $w$ we have
$$w(\gamma_0')=2\int_{\vc{u}_1}^B g(z)\; dz\geq \int_A^B g(z)\; dz.$$
\end{pfof}

To prove Proposition \ref{vsplit} we need Lemmas \ref{vshift}, \ref{vmono} below.

\begin{lemma}
\label{vshift}
If $\vc{a}\leq\vc{b}$ and $\vc{k}\in[0,\infty)^2$ then
$W(\vc{a},\vc{b})\geq W(\vc{a}+\vc{k},\vc{b}+\vc{k})$.
\end{lemma}
(Note that here $\vc{k}$ might not be an element of $\rtp$).

\begin{pf}
We have
\begin{eqnarray*}
W(\vc{a}+\vc{k},\vc{b}+\vc{k})
&=&\inf_{\gamma:\vc{a}+\vc{k}\rightarrow\vc{b}+\vc{k}} w(\gamma) \\
&=&\inf_{\gamma:\vc{a}\rightarrow\vc{b}} w(\gamma+\vc{k}) \\
&\leq&\inf_{\gamma:\vc{a}\rightarrow\vc{b}} w(\gamma) \qquad\text{ since g is decreasing}\\
&=&W(\vc{a},\vc{b}),
\end{eqnarray*}
where $\gamma+\vc{k}$ denotes the shifted path obtained by adding $\vc{k}$ to each point in $\gamma$.
\end{pf}

\begin{lemma}
\label{vmono}
If $\vc{a}\leq\vc{b}\leq\vc{c}$ then
$W(\vc{a},\vc{c})\geq W(\vc{b},\vc{c})$.
\end{lemma}
We prove Lemma \ref{vmono} via the following.
\begin{lemma}
\label{vmono2}
If $\vc{a}\leq\vc{c}$ and either $\vc{b}=(\vc{a}_1,\vc{c}_2)$ or $\vc{b}=(\vc{c}_1,\vc{a}_2)$ then
$W(\vc{a},\vc{c})\geq W(\vc{b},\vc{c})$.
\end{lemma}

\begin{pf}
Without loss of generality suppose that $\vc{b}=(\vc{a}_1,\vc{c}_2)$.  Let $\gamma$ be any path from $\vc{a}$ to $\vc{b}$, and let $\delta$ be the unique (straight, horizontal) path from $\vc{b}$ to $\vc{c}$.  Since $g$ is decreasing and positive we have
$$W(\vc{b},\vc{c})=w(\delta)=\int_\delta g(y)dx\leq \int_\gamma g(y)dx\leq w(\gamma).$$
\end{pf}

\begin{pfof}{Lemma \ref{vmono}}
Let $\gamma$ be any path from $\vc{a}$ to $\vc{b}$.  Without loss of generality suppose that $\{\vc{b}\}\succeq\gamma$.  Let $\vc{u}=(\vc{a}_1,\vc{b}_2)$ and $\vc{v}=(\vc{v}_1,\vc{b}_2)$ where
$$\vc{v}_1=\inf\{t:(t,\vc{b}_2)\in\gamma\}.$$
Write $\gamma=\gamma_1\cup\gamma_2$ where $\gamma_1$ is a path $\vc{a}$ to $\vc{v}$ and $\gamma_2$ is a path $\vc{v}$ to $\vc{b}$ (both of which are thus uniquely defined).  Let $\delta=\delta_1\cup\delta_2$ where $\delta_1$ is the unique path from $\vc{u}$ to $\vc{b}$ and $\delta_2$ is the unique path from $\vc{b}$ to $\vc{v}$.  See Figure \ref{monopf}.
\begin{figure}
\centering
\begin{picture}(0,0)%
\includegraphics{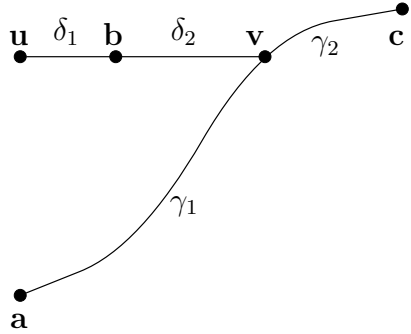}%
\end{picture}%
\setlength{\unitlength}{0.00041700in}%
\begin{picture}(4966,4085)(1118,-4363)
\put(1201,-4336){\makebox(0,0)[b]{$\vc{a}$}}
\put(6001,-736){\makebox(0,0)[b]{$\vc{c}$}}
\put(2401,-736){\makebox(0,0)[b]{$\vc{b}$}}
\put(4201,-736){\makebox(0,0)[b]{$\vc{v}$}}
\put(1201,-736){\makebox(0,0)[b]{$\vc{u}$}}
\put(3301,-2911){\makebox(0,0)[b]{$\gamma_1$}}
\put(5101,-886){\makebox(0,0)[b]{$\gamma_2$}}
\put(1801,-736){\makebox(0,0)[b]{$\delta_1$}}
\put(3301,-736){\makebox(0,0)[b]{$\delta_2$}}
\end{picture}
\caption{Illustration of the proof of Lemma \ref{vmono}.}
\label{monopf}
\end{figure}

Now
\begin{eqnarray*}
w(\gamma)&=&w(\gamma_1)+w(\gamma_2) \\
&\geq& w(\delta)+w(\gamma_2) \quad\text{ by Lemma \ref{vmono2}}\\
&\geq& w(\delta_2)+w(\gamma_2) \\
&=&w(\delta_2\cup\gamma_2) \\
&\geq &W(\vc{b},\vc{c}) \quad\text{ since $\delta_2\cup\gamma_2$ is a path from $\vc{b}$ to $\vc{c}$.}
\end{eqnarray*}
\end{pfof}

We are now ready to prove Proposition \ref{vsplit}; we start with a version without the ``error term''.

\begin{lemma}
\label{split2}
Suppose that $\vc{a}\leq\vc{b}$; $\vc{c}\leq\vc{d}$; $\vc{r}\geq\vc{b}$; $\vc{r}\geq\vc{d}$ and $\vc{r}\leq\vc{b}+\vc{d}$.  Then there exists $\vc{s}$ satisfying $\vc{s}\leq \vc{r}$ and $\vc{s}\leq \vc{a}+\vc{c}$ such that
$$W(\vc{a},\vc{b})+W(\vc{c},\vc{d})\geq W(\vc{s},\vc{r}).$$
Indeed we may take $\vc{s}=\vc{a}\vee(\vc{a}-\vc{b}+\vc{c}-\vc{d}+\vc{r})$.
\end{lemma}

\begin{pf}
Define $\vc{s}$ as indicated, and note that $\vc{s}\geq\vc{a}$; $\vc{s}\leq\vc{a}+\vc{c}$; and $\vc{s}\leq \vc{r}$.  Then we have
\begin{eqnarray*}
\lefteqn{W(\vc{a},\vc{b})+W(\vc{c},\vc{d})} \\
&\geq &W(\vc{s},\vc{s}-\vc{a}+\vc{b})+W(\vc{r}+\vc{c}-\vc{d},\vc{r}) \quad\text{by Lemma \ref{vshift} twice} \\
&\geq &W(\vc{s},\vc{s}-\vc{a}+\vc{b})+W(\vc{s}-\vc{a}+\vc{b},\vc{r}) \quad\text{by Lemma \ref{vmono}} \\
&\geq & W(\vc{s},\vc{r}) \quad\text{by Lemma \ref{vsum}}.
\end{eqnarray*}
In the first inequality we have used the facts that $\vc{s}-\vc{a}\geq (0,0)$ and $\vc{r}-\vc{d}\geq (0,0)$.  For the second (and third) we must check that $\vc{r}+\vc{c}-\vc{d}\leq \vc{s}-\vc{a}+\vc{b}\leq \vc{r}$, which is achieved as follows:
$$\vc{r}-(\vc{s}-\vc{a}+\vc{b})=(\vc{r}+\vc{a}-\vc{b})-\vc{s}=(\vc{r}-\vc{b})\wedge (\vc{d}-\vc{c})\geq (0,0),$$
and
$$(\vc{s}-\vc{a}+\vc{b})-(\vc{r}+\vc{c}-\vc{d})=\vc{s}-(\vc{a}-\vc{b}+\vc{c}-\vc{d}+\vc{r})=(0,0)\vee(\vc{b}-\vc{c}+\vc{d}-\vc{r})\geq (0,0).$$
\end{pf}

\begin{pfof}{Proposition \ref{vsplit}}
Let $\vc{r}'=\vc{r}\wedge (\vc{b}+\vc{d})$, and note that $\vc{r'}\geq\vc{b}$; $\vc{r'}\geq\vc{d}$; $\vc{r'}\leq\vc{b}+\vc{d}$; $\vc{r'}\leq\vc{r}$ and $\vc{r'}\geq\vc{r}-(q,q)\geq (2Z-q,2Z-q)\geq (Z,Z)$.  Hence we may apply Lemma \ref{split2} to obtain
\begin{eqnarray*}
W(\vc{a},\vc{b})+W(\vc{c},\vc{d})&\geq &W(\vc{s},\vc{r}') \\
&\geq &W(\vc{s},\vc{r})-W(\vc{r}',\vc{r}) \quad\text{ by Proposition \ref{vsum}} \\
&\geq &W(\vc{s},\vc{r})-2qg(Z) \quad\text{ by Proposition \ref{vupper}},
\end{eqnarray*}
where $\vc{s}=\vc{a}\vee(\vc{a}-\vc{b}+\vc{c}-\vc{d}+\vc{r'})=\vc{a}\vee[(\vc{a}+\vc{c})\wedge (\vc{a}+\vc{c}+\vc{r}-\vc{b}-\vc{d})]$, and $\vc{s}\leq\vc{r}'\leq\vc{r}$ and $\vc{s}\leq \vc{a}+\vc{c}$.
\end{pfof}

\section{Border Events}
\label{secbord}

\begin{figure}
\centering
\begin{picture}(0,0)%
\includegraphics{border.pstex}%
\end{picture}%
\setlength{\unitlength}{0.00041700in}%
\begin{picture}(5702,3644)(579,-3383)
\put(3301,-1636){\makebox(0,0)[b]{$R$}}
\put(5251,-1636){\makebox(0,0)[b]{$R_4$}}
\put(5251,-286){\makebox(0,0)[b]{$R_5$}}
\put(3301,-286){\makebox(0,0)[b]{$R_6$}}
\put(1201,-1636){\makebox(0,0)[b]{$R_8$}}
\put(1201,-286){\makebox(0,0)[b]{$R_7$}}
\put(6151,-1636){\makebox(0,0)[b]{$R'$}}
\put(3301,-3136){\makebox(0,0)[b]{$R_2$}}
\put(1201,-3136){\makebox(0,0)[b]{$R_1$}}
\put(5251,-3136){\makebox(0,0)[b]{$R_3$}}
\end{picture}
\caption{The rectangles $R_1,\ldots, R_8$.}
\label{r1r8}
\end{figure}
Let $R,R'$ be two rectangles satisfying $R\subseteq R'$.  Define rectangles $R_1,\ldots ,R_8$ (some of which may be empty) according to Figure \ref{r1r8}, so that $R'$ is the disjoint union of $R,R_1,\ldots,R_8$.  Define $D(R,R')$ to be the event that each of the two rectangles $R_1\cup R_8 \cup R_7$ and $R_3\cup R_4\cup R_5$ is horizontally traversable, and each of the two rectangles $R_1\cup R_2 \cup R_3$ and $R_7\cup R_6\cup R_5$ is vertically traversable.  One may think of $D(R,R')$ as a necessary condition for the ``border'' between $R$ and $R'$ to be ``traversable'' from $R$ to $R'$.
More precisely, the event has the following properties.
\begin{mylist}
\item
If $R'$ is internally spanned then $D(R,R')$ occurs.
\item
$D(R,R')$ is defined in terms of the states of sites in $R'\setminus R$.
\end{mylist}
Property (i) holds because if $R'$ is internally spanned then it must be horizontally and vertically traversable.
The purpose of this section is to prove the following upper bound on the probability of $D(R,R')$.

\begin{prop}
\label{border}
For any $Z>0$ and $c\in(0,1/2)$, there exist $Q=Q(c,Z)<\infty$ and $T=T(c,Z)\in(0,Z/2)$ such that for any rectangles $R\subseteq R'$ with dimensions $(m,n)$ and $(m+s,n+t)$ respectively, and any $q>0$, provided $m,n\geq Z/q$ and $s,t\leq T/q$ we have
$$P_p(D(R,R'))\leq Q \exp \bigg(-(1-2c)\left[g(nq)s+g(mq)t\right]\bigg).$$
\end{prop}
In the applications of Proposition \ref{border}, it will be essential that $Q,T$ do not depend on $q$.

\begin{pf}
Let $H$ be the ``corner region'', $H=R_1\cup R_3\cup R_5\cup R_7$, and let $Y$ be the number of occupied sites in $H$, $Y=|X\cap H|$.  The idea of the proof is as follows: if the four events in the definition of $D$ were independent, the proof would be easy.  If $Y$ is small, then the events are nearly independent, and if $s,t$ are sufficiently small compared with $m,n$, the probability $Y$ is large can be made smaller than the bound we are trying to obtain for $P_p(D)$.

Without loss of generality suppose that $s\leq t$.  We have
\begin{eqnarray}
\lefteqn{P(D)} \nonumber\\
&=&P(D\mid Y \leq cs)P(Y\leq cs) + P(D\mid cs<Y \leq ct)P(cs<Y\leq ct) \nonumber\\
&&\quad + P(D \mid Y>ct)P(Y>ct) \nonumber\\
&\leq& P(D\mid Y \leq cs) + P(D\mid cs<Y \leq ct)P(Y>cs) \label{split} +P(Y>ct)
\end{eqnarray}
We claim that the terms appearing on the right side of (\ref{split}) may be bounded as follows
\begin{eqnarray}
P(D\mid Y \leq cs)&\leq& e^{4g(Z)}\exp -(1-2c)[g(nq)s+g(mq)t] \label{d1}\\
P(D\mid cs<Y \leq ct)&\leq& e^{2g(Z)}\exp -(1-2c)g(mq)t \label{d2} \\
P(Y>cs)&\leq&\exp - cs(\log c -\log T -1) \label{x1} \\
P(Y>ct)&\leq&\exp - ct(\log c -\log T -1) \label{x2}
\end{eqnarray}
for a suitable choice of $T$.

The inequality (\ref{x1}) is an instance of the Chernoff bound, as follows.  $Y$ is a binomial random variable with parameters $st$ and $p$, so
\begin{eqnarray*}
P(Y>cs)&=&P(e^{a(Y-cs)}>1) \quad\text{ for any }a>0\\
&\leq&E(e^{a(Y-cs)})\\
&=&\exp[-csa +st\log(1-p+pe^a)]\\
&\leq&\exp[-csa+stpe^a]\\
&=&\exp[-cs\log(c/(tp))+cs] \quad\text{ taking }e^a=c/(tp)\\
&\leq&\exp-cs(\log c-\log T-1),
\end{eqnarray*}
provided $tp\leq tq\leq T$.  ($T$ will be a function of $Z$ and $c$ to be chosen later).  The bound (\ref{x2}) follows similarly.

To prove (\ref{d2}), we condition further on which sites in $H$ are occupied.  If for each choice of this set of sites (satisfying $cs<Y \leq ct$) the conditional probability of $D$ is bounded above by the right side of (\ref{d2}), the desired bound will follow.  For simplicity, suppose first that $Y=ct$, and also that the $ct$ occupied sites in $H$ all lie in different horizontal rows.  These sites (together with $R$) split $R_2\cup R_6$ into $ct+2$ ``horizontal strips'', some of which may be empty.  (More precisely, if we remove from $R_2\cup R_6$ every row which contains an occupied site in $H$, we are left with $ct+2$ rectangles of width $m$, some of which may be empty).  In order for $D$ to occur, a necessary condition is that each of these strips is vertically traversable.  By Lemma \ref{traversable} (i), the probability of this is at most
$$\exp -[t-ct-(ct+2)]g(mq)$$
(since the sum of the vertical heights of the strips is $t-ct$, hence the sum of their heights minus one is $t-ct-(ct+2)$).  Provided $mq\geq Z$, the above expression is at most
\begin{equation}
e^{2g(Z)}\exp -(1-2c)g(mq)t,
\label{horbound}
\end{equation}
since $g$ is decreasing.
Now, if $Y<ct$, or if some of the occupied sites in $H$ lie in the same horizontal rows, clearly the conditional probability of $D$ will be even smaller, hence we have proved (\ref{d2}).

The bound (\ref{d1}) is proved similarly.  If $Y\leq cs$ then $Y\leq ct$ also, and by conditioning on the occupied sites in $H$ we obtain a collection of horizontal strips in $R_2\cup R_6$ together with a collection of vertical strips in $R_4\cup R_8$.  For $D$ to occur, each of horizontal strips must be vertically traversable, and each of vertical strips must be horizontally traversable, but these two events are independent (conditional on the set of occupied sites in $H$), so the two bounds corresponding to (\ref{horbound}) are multiplied, to obtain the right side of (\ref{d1}).

Now choose $T>0$ sufficiently small that
$$c(\log c-\log T-1)\geq 2(1-2c)g(Z)$$
(and also $T<Z/2$).
Since $g$ is decreasing this ensures that
\begin{equation}
\exp - cs(\log c -\log T -1) \leq \exp -(1-2c)g(nq)s
\label{comp1}
\end{equation}
and also
\begin{equation}
\exp - ct(\log c -\log T -1) \leq \exp -(1-2c)[g(nq)s+g(mq)t].
\label{comp2}
\end{equation}
In the latter we have used the fact that $s\leq t$ so $2g(Z)t \geq g(Z)s+g(Z)t$.  

Now substituting (\ref{d1})--(\ref{x2}) into (\ref{split}) and using (\ref{comp1}),(\ref{comp2}) we obtain
$$P(D)\leq 3e^{4g(Z)} \exp -(1-2c)[g(nq)s+g(mq)t].$$
\end{pf}

\section{Disjoint Spanning}
\label{secdisj}

For a collection of increasing events $A_1,\ldots ,A_k$ on $\{0,1\}^{\zt}$, the event $A_1\circ\cdots\circ A_k$ that $A_1,\ldots ,A_k$ {\dof occur disjointly} is defined as the event that there exist pairwise disjoint full sets of sites $K_1,\ldots ,K_k$ such that for each $i$, $A_i$ occurs whenever $K_i$ is full.  The BK inequality states that if $A_1,\ldots ,A_k$ are defined in terms the states of a finite set of sites then
$$P_p(A_1\circ\cdots\circ A_k)\leq P_p(A_1)\cdots P_p(A_k).$$
For more details see for example \cite{g2} p37.

\begin{figure}
\centering \resizebox{!}{1.5in}{\includegraphics{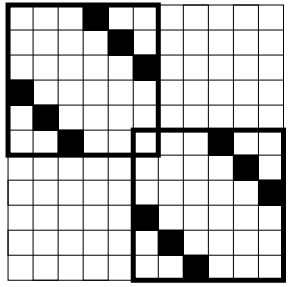}}
\caption{In this example the rectangles $R',R''$ must overlap.}
\label{overlap}
\end{figure}

\begin{prop}
\label{disjoint}
Let $R$ be a rectangle with $|R|\geq 2$.
If $R$ is internally spanned then there exist two distinct non-empty rectangles $R',R''$ such that
\begin{mylist}
\item
the {\em strict} inclusions $R'\subset R$, $R''\subset R$ hold,
\item
$\langle R'\cup R''\rangle =R$,
\item
$\{R' \text { is internally spanned}\}\circ\{R'' \text { is internally spanned}\}$ occurs.
\end{mylist}
\end{prop}
To see the subtlety of Proposition \ref{disjoint}, note that we cannot in general take the two rectangles $R',R''$ to be disjoint; in Figure \ref{overlap} for example, the whole square $R$ is internally spanned, but the only possible choice for $R',R''$ is the pair of 6 by 6 squares indicated.  It is for this reason that the concept of disjoint occurrence is important.
The idea of the proof of Proposition \ref{disjoint} is simple: we run an algorithm which produces successively larger internally spanned rectangles, combining two rectangles into one at each step and eventually obtaining $R$; then we consider the last step.  The details of the proof require a little more care.

\begin{pfof}{Proposition \ref{disjoint}}
If $K$ is any finite set of sites, we may construct $\langle K\rangle$ via the following algorithm.  For each time step $t=0,1,\ldots ,\tau$, we shall construct a collection of $m_t$ rectangles $R_1^t,\ldots , R_{m_t}^t$, and corresponding sets of sites $K_1^t,\ldots , K_{m_t}^t$, with the following properties:
\begin{mylist}
\item
$K_1^t,\ldots , K_{m_t}^t$ are pairwise disjoint;
\item
$K_i^t\subseteq K$;
\item
$R_i^t=\langle K_i^t\rangle$;
\item
if $i\neq j$ then $R_i^t \not\subseteq R_j^t$;
\item
$K\subseteq H^t \subseteq \langle K\rangle$, where
$$H^t:=\bigcup_{i=1}^{m_t} R_i^t.$$
\end{mylist}

Initially, the rectangles and sets of sites are just the individual sites of $K$.  That is, let $K$ be enumerated as $K=\{x_1,\ldots, x_k\}$, and set 
$m_0=k$ \text{ and } $R_i^0=K_i^0=\{x_i\},$
so that in particular
$$H^0 =K.$$
The final set of rectangles $R_1^\tau,\ldots R_{m_\tau}^\tau$ will have the property that
\begin{equation}
\label{algfin}
H^\tau =\langle K\rangle.
\end{equation}

Before describing the algorithm we make the following observation.  Let $R_1,R_2$ be two distinct rectangles neither of which is a subset of the other, and consider $\langle R_1\cup R_2\rangle$.  The following three possibilities exist.
\begin{mylist}
\item[(a)]
$\langle R_1\cup R_2\rangle=R_1\cup R_2$, and $R_1\cup R_2$ is not connected.
\item[(b)]
$\langle R_1\cup R_2\rangle=R_1\cup R_2$, and $R_1\cup R_2$ is a rectangle.
\item[(c)]
$\langle R_1\cup R_2\rangle\supset R_1\cup R_2$, and $R_1\cup R_2$ is a rectangle.
\end{mylist}
(As usual, a set of sites is said to be connected if any two sites can be joined via a sequence of sites at Euclidean distances 1).

The algorithm proceeds as follows.  Suppose $R_1^t,\ldots , R_{m_t}^t$ and $K_1^t,\ldots , K_{m_t}^t$ have already been constructed.
\begin{mylist}
\item[Step (I).]  If there do not exist a pair of rectangles $R_i^t,R_j^t$ with $i\neq j$ such that $\langle R_i^t \cup R_j^t\rangle$ is a rectangle (that is, if case (a) above holds for all pairs), then {\bf STOP}, and set $\tau=t$.

\item[Step (II).]  If there do exist $R_i^t,R_j^t$ with $i\neq j$ such that $\langle R_i^t \cup R_j^t\rangle$ is a rectangle (case (b) or (c) above), then choose one such pair of rectangles.  Denote the rectangle $\langle R_i^t \cup R_j^t\rangle$ by $R'$.  Also let $K'=K_i^t \cup K_j^t$.

\item[Step (III).]  Construct the state $(R_1^{t+1},K_1^{t+1}),\ldots ,(R_{m_{t+1}}^{t+1},K_{m_{t+1}}^{t+1})$ at time $t+1$ as follows.  From the list $(R_1^{t},K_1^{t}),\ldots ,(R_{m_{t}}^{t},K_{m_{t}}^{t})$ at time $t$, delete every pair $(R_l^t,K_l^t)$ for which $R_l^t\subseteq R'$.  This includes $(R_i^t,K_i^t)$ and $(R_j^t,K_j^t)$, and may include others.  Then add $(R',K')$ to the list.
\item[Step (IV).] Increase $t$ by 1 and return to Step (I).
\end{mylist}

It is straightforward to see that properties (i)--(v) are preserved by this procedure.  Also $m_t$ is strictly decreasing with $t$, so the algorithm must stop eventually.  To check that (\ref{algfin}) is satisfied, observe that if there exists a site $x\in \langle K\rangle \setminus H^\tau$, then there must exist a site $y\in \langle K\rangle \setminus H^\tau$ having at least two neighbours in $H^\tau$, but these neighbours must lie in distinct rectangles $R_i^\tau,R_j^\tau$, and hence the algorithm should not have stopped at time $\tau$. 

Furthermore, observe that if $\langle K\rangle$ is a single rectangle $R$, then we must have $m_\tau=1$ and $R_1^\tau=R$.  If not, since $H^\tau=R$, there must exist two distinct rectangles $R_i^\tau,R_j^\tau$ whose union is connected, and again this means that the algorithm should not have stopped.

Finally, to prove the proposition, note that if $R$ is internally spanned then running the algorithm on the set of sites $K=R\cap X$ results in $m_\tau=1$ and $R_1^\tau=R$.  If $|R|\geq 2$ then there must have been at least one step, $\tau\geq 1$.  Now considering the last time step of the algorithm (from time $\tau-1$ to time $\tau$) we obtain the two rectangles $R'=R_i^{\tau-1}$, $R''=R_j^{\tau-1}$ with all the required properties.
\end{pfof}

\section{Hierarchies}
\label{sechier}

A {\dof directed graph} is a set of {\dof vertices} $V$ together with a set of ordered pairs of vertices $E$.  If $(u,v)\in E$ then we say there is an {\dof edge} from $u$ to $v$ and write $u\leadsto v$.  The {\dof children} of a vertex $u$ are all the vertices $v$ such that $u\leadsto v$.

\begin{figure}
\centering
\begin{picture}(0,0)%
\includegraphics{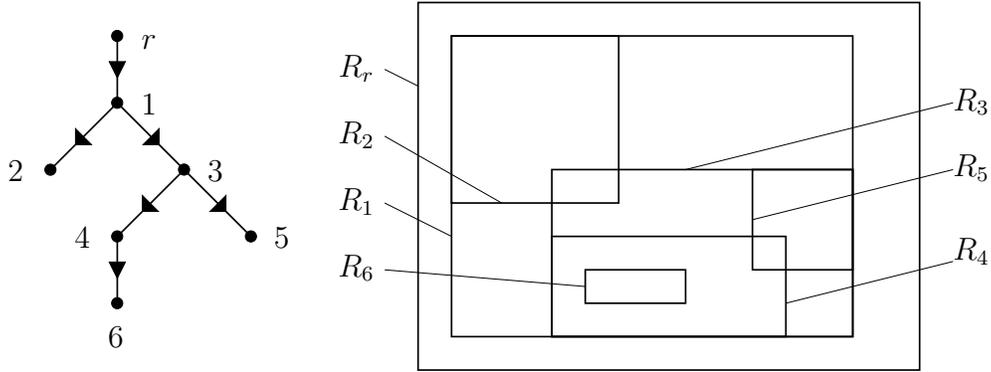}%
\end{picture}%
\setlength{\unitlength}{0.00058300in}%
\begin{picture}(8855,3344)(532,-2783)
\put(9226,-1786){\makebox(0,0)[b]{$R_4$}}
\put(9226,-1036){\makebox(0,0)[b]{$R_5$}}
\put(9226,-436){\makebox(0,0)[b]{$R_3$}}
\put(3676,-1936){\makebox(0,0)[b]{$R_6$}}
\put(3676,-1336){\makebox(0,0)[b]{$R_1$}}
\put(3676,-736){\makebox(0,0)[b]{$R_2$}}
\put(3676,-136){\makebox(0,0)[b]{$R_r$}}
\put(1802,163){\makebox(0,0)[b]{$r$}}
\put(1802,-437){\makebox(0,0)[b]{$1$}}
\put(602,-1037){\makebox(0,0)[b]{$2$}}
\put(2402,-1037){\makebox(0,0)[b]{$3$}}
\put(1202,-1637){\makebox(0,0)[b]{$4$}}
\put(3002,-1637){\makebox(0,0)[b]{$5$}}
\put(1502,-2537){\makebox(0,0)[b]{$6$}}
\end{picture}
\caption{An example of a hierarchy.}
\label{exhier}
\end{figure}

A {\dof hierarchy} ${\cal H}$ is a finite directed graph in which every vertex $v$ is labeled with a non-empty rectangle $R_v$ (where the rectangles corresponding to different vertices are not necessarily distinct), with the following properties.  The graph is a tree.  There is a special vertex $r$ called the {\dof root}, and all edges are directed away from the root. (So for any vertex $v$ there is a unique directed chain of vertices $r=u_0\leadsto u_1\leadsto \cdots \leadsto u_k=v$).  If $u\leadsto v$ then we have the {\em strict} inclusion $R_u\supset R_v$. 
Every vertex has $0$, $1$ or $2$ children.  A vertex with no children is called a {\dof seed}.  If $u$ has exactly one child $v$, we call $u$ {\dof normal} and write $u\Rightarrow v$.  If $u$ has two children $v,w$ we call $u$ a {\dof splitter} and write $u\rightrightarrows (v,w)$.  If $u\rightrightarrows (v,w)$ then we have $\langle R_v \cup R_w\rangle = R_u$.
See Figure \ref{exhier} for an example of a hierarchy.

\pagebreak
We say that the hierarchy ${\cal H}$ {\dof occurs} if all the following events occur {\em disjointly}:
$$D(R_v,R_u) \qquad\text{ for each pair } u,v \text{ such that }u\Rightarrow v,$$
\hspace{\parindent} and
$$\{R_w \text{ is internally spanned}\} \qquad\text{ for each seed } w.$$
(Note in particular that if for example $R_v=R_w$ for two distinct seeds $v,w$ then we require that the event $\{R_v \text{ is internally spanned}\}$ occurs twice disjointly; that is, there are two disjoint sets of occupied sites each of which internally spans $R_v$).

Fix $q,T,Z$ satisfying $0<4q\leq 2T\leq Z \leq 1/2$ (in fact we shall be concerned with the case $0<q\ll T\ll Z\ll 1$).  We call a hierarchy ${\cal H}$ {\dof good} if it satisfies all the following.
\begin{mylist}
\item
If $w$ is a seed then $$\sht(R_w)<2Z/q;$$
\item
if $u$ is normal or a splitter then $$\sht(R_u)\geq 2Z/q;$$
\item
if $u\Rightarrow v$ and $v$ is a seed or a normal vertex then $$\phi(R_u)-\phi(R_v)\in[T/(2q),T/q];$$
\item
if $u\Rightarrow v$ and $v$ is a splitter then $$\phi(R_u)-\phi(R_v)\leq T/q;$$
\item
if $u\rightrightarrows(v,w)$ then 
\begin{eqnarray*}
\phi(R_u)-\phi(R_v)&\geq &T/(2q)\\
 \text{and }\quad \phi(R_u)-\phi(R_w)& \geq & T/(2q).
\end{eqnarray*}
\end{mylist}

\begin{prop}
\label{hier}
Let $q,T,Z$ satisfy $0<4q\leq 2T\leq Z\leq 1/2$, and let $R$ be a non-empty rectangle.  If $R$ is internally spanned then some good hierarchy with root-label $R_r=R$ occurs.
\end{prop}

\begin{pf}
The proof is by induction on the size of the rectangle.  Let $R$ be a rectangle and suppose the proposition holds for all rectangles with semi-perimeter less than $\phi(R)$.

If $\sht(R)< 2Z/q$, then the good hierarchy having only one vertex $r$ (which is the root and a seed), and $R_r=R$, occurs.

If $\sht(R)\geq 2Z/q$ then we construct a sequence of non-empty internally spanned rectangles $R=S_0\supset S_1\supset \cdots\supset S_m$ by an algorithm as follows.  The idea is that $S_1,\ldots ,S_m$ are successive attempts to find a rectangle $S$ such that $\phi(R)-\phi(S)$ is in the range $[T/(2q),T/q]$; the attempt may succeed in which case the root will be a normal vertex, or we may ``overshoot'' in which case we need to introduce a splitter.  Here are the details.  Given $S_i$, apply Proposition \ref{disjoint} to obtain two rectangles $S_i',S_i''$.  Let $S_{i+1}$ be the one of $S_i',S_i''$ with the larger semi-perimeter $\phi$  (choosing according to an arbitrary rule if they are equal).  Stop, after $m\geq 1$ steps, the first time $\phi(R)-\phi(S_m)\geq T/(2q)$.  (This must occur eventually because the sequence of rectangles is strictly decreasing, and $\phi(R)\geq 4Z/q > T/(2q)+2$, so that a rectangle $S$ containing only one site does satisfy $\phi(R)-\phi(S)\geq T/(2q)$; and we must have $m\geq 1$ because $T/(2q)>0$).

Now consider the following three possible cases.
\begin{mylist}
\item
If $\phi(R)-\phi(S_m)\leq T/q$, \\[\parsep]
 then we have that $\phi(R)-\phi(S_m)\in [T/(2q),T/q]$, $S_m\subset R$, and $R$ and $S_m$ are internally spanned.  By the inductive hypothesis, there exists a good hierarchy ${\cal H}'$ with root $r'$ and root-label $R_{r'}=S_m$.  Furthermore, the events $D(S_m,R)$ and $\{{\cal H}' \text{ occurs}\}$ occur disjointly, since they are defined in terms of disjoint sets of sites.  We construct a hierarchy ${\cal H}$ as follows: start with ${\cal H}'$, and add a new vertex $r$, with $R_r=R$, and an edge from $r$ to $r'$ so that $r\Rightarrow r'$.  It follows from the above observations that ${\cal H}$ is good and occurs.
\item
If $\phi(R)-\phi(R_m)> T/q$ and $m=1$, \\[\parsep]
 then there exist disjointly internally spanned rectangles $S_0',S_0''$ such that $\langle S_0' \cup S_0''\rangle=R$, and $\phi(R)-\phi(S_0'), \phi(R)-\phi(S_0'')>T/q>T/(2q)$.  By Proposition \ref{disjoint} and the inductive hypothesis there exist disjointly occuring good hierarchies ${\cal H}',{\cal H}''$ with roots $r',r''$ and root-labels $R_{r'}=S_0', R_{r''}=S_0''$, and we may take the vertex sets of ${\cal H}',{\cal H}''$ to be disjoint.  Now we construct ${\cal H}$ from ${\cal H}'\cup {\cal H}''$  (this last object being a labeled directed graph defined in the obvious way) by adding a new vertex $r$, with $R_r=R$, and edges so that $r\rightrightarrows (r',r'')$.  It is easily seen that ${\cal H}$ is a good hierarchy and occurs.
\item
If $\phi(R)-\phi(R_m)> T/q$ and $m\geq 2$, \\[\parsep]
 then we have internally spanned rectangles $S_{m-1},S_{m-1}',S_{m-1}''$ satisfying $R\supset S_{m-1}$ and $S_{m-1}=\langle S_{m-1}' \cup S_{m-1}''\rangle$, where $\phi(R)-\phi(S_{m-1})<T/(2q)$ and $\phi(R)-\phi(S_{m-1}'),\phi(R)-\phi(S_{m-1}'')>T/q$, and therefore we have $\phi(S_{m-1})-\phi(S_{m-1}'),\phi(S_{m-1})-\phi(S_{m-1}'')\geq T/(2q)$.  By Proposition \ref{disjoint} and the inductive hypothesis there exist disjointly occuring good hierarchies ${\cal H}',{\cal H}''$ with roots $r',r''$ and root-labels $R_{r'}=S_{m-1}', R_{r''}=S_{m-1}''$, and we may take the vertex sets of ${\cal H}',{\cal H}''$ to be disjoint.  Now we construct ${\cal H}$ from ${\cal H}'\cup {\cal H}''$ by adding new vertices $r,y$ with $R_r=R, R_y=S_{m-1}$, and new edges so that $r\Rightarrow y$ and $y\rightrightarrows (r',r'')$.  Then it is easily seen that ${\cal H}$ is a good hierarchy and occurs.
\end{mylist}
\end{pf}

\section{Upper Bound}
\label{secupper}

We are now ready to prove Theorem \ref{mainprop} (ii).
Fix $B>2$, and let $A=c=1/B$.  We shall prove that
$$\liminf_{p\rightarrow 0} -p\log I(\lfloor B/p\rfloor,p) \geq 2(1-2c)\int_A^B g(z)\; dz,$$
from which Theorem \ref{mainprop} (ii) follows.  The approach is to use Proposition \ref{hier}, and obtain upper bounds on the number of possible good hierarchies, and on the probability that each one occurs.

Choose $Z>0$ sufficiently small that $Z< A/2$ and
\begin{equation}
\label{zchoice}
g(2Z)\geq \frac{4\lambda}{A}
\end{equation}
where $\lambda=\pi^2/18$.
Recall that $g(z)\rightarrow\infty$ as $z\rightarrow 0$, so this is indeed possible; the reason for this particular choice of $Z$ will become clear later.
Finally choose $T=T(c,Z),Q=Q(c,Z)$ according to Proposition \ref{border}.  It will also be convenient to assume that $q<T/2$.  Thus we have
$$16q<8T<4Z<2A<1<B/2.$$
We shall be concerned with ``large'' $B$ and ``very small'' $q$, in which case we have
$$q\ll T\ll Z \ll A \ll 1 \ll B.$$
Later we will let $q\rightarrow 0$ while keeping $B$ fixed.  It will be important to distinguish between quantities which depend only on $B$ (such as $A,c,Z,T,Q$) and those which also depend on $q$.

If $R$ is a rectangle we define 
$$V(R)=q\lng(R)\; g\bigg(q\sht(R)\bigg),$$
and for rectangles $R\subseteq R'$ we define
$$U(R,R')=W\bigg(q\dim(R),q\dim(R')\bigg),$$
where $W$ is defined as in Section \ref{secvar}.
We claim that for any $R$ we have
\begin{equation}
\label{energyv}
P_p(R \text{ is internally spanned})\leq \exp-q^{-1} V(R),
\end{equation}
while for any $R\subseteq R'$ such that $\sht(R') \geq 2Z/q$ and $\sem(R')-\sem(R)\leq T/q$ we have
\begin{equation}
\label{energyu}
P_p(D(R,R'))\leq Q \exp -q^{-1} (1-2c)U(R,R').
\end{equation}

Inequality (\ref{energyv}) follows from Lemma \ref{traversable} (ii) and Lemma \ref{travapps} (i) (using which\-ever of East- and North- traversability gives the better bound).  For (\ref{energyu}), note that $\dim(R')-\dim(R) \leq (T/q,T/q)$, and $\sht(R)\geq 2Z/q- T/q\geq Z/q$, so $R,R'$ satisfy the conditions of Proposition \ref{border}.  Combining this with Proposition \ref{vupper} and the definition of $U$ above yields (\ref{energyu}).

Now if ${\cal H}$ is any good hierarchy, by BK inequality and the definition of ${\cal H}$ occuring, we have  
\begin{eqnarray}
\lefteqn{P_p({\cal H} \text{ occurs})} \nonumber \\
&\leq& \prod_{u\Rightarrow v} P_p\bigg(D(R_v,R_u)\bigg) \prod_{w \text{ seed}} P_p(R_w \text{ internally spanned}) \nonumber \\
&\leq& Q^{\nnorm} \exp -q^{-1}
\left[(1-2c)\sum_{u\Rightarrow v} U(R_v,R_u)+\sum_{w \text{ seed}} V(R_w)\right], \label{hocc}
\end{eqnarray}
where the first product and sum are over all pairs of vertices $u,v$ such that $u\Rightarrow v$, the second product and sum are over all seeds $w$, and 
$\nnorm$ is the number of normal vertices of ${\cal H}$.
The second inequality in (\ref{hocc}) follows from (\ref{energyv}),(\ref{energyu}) above and from properties (ii),(iii),(iv) of a good hierarchy.

Next we derive lower bounds for the two sums in (\ref{hocc}). 
\begin{lemma}
\label{pod}
For any good hierarchy ${\cal H}$ with root-label $R_r=R$, there exists a rectangle $S=S({\cal H})\subseteq R$ satisfying
$$\dim(S)\leq \sum_{w \text{{\rm\ seed}}} \dim(R_w)$$
such that
$$\sum_{u\Rightarrow v} U(R_v,R_u) \geq U(S,R) - \nsplit 2qg(Z),$$
where $\nsplit$ is the number of splitters of ${\cal H}$.
\end{lemma}
We call $S$ as above the {\dof pod} of ${\cal H}$.  The idea is that, if the total size of the seeds is not too big, then the pod is not too big, and therefore the first sum in (\ref{hocc}) is large enough to give a good bound.  Note that the location of the pod is actually immaterial - only its dimensions are ever used.  It is defined to be a rectangle rather than just a 2-vector as a notational convenience only.  The pod is function of the hierarchy ${\cal H}$ only, not of $B$ or $q$ (even though the definitions of a good hierarchy and of $U(\cdot,\cdot)$ do depend on $B$ and $q$).

The following will be needed in the proof of Lemma \ref{pod}.
\begin{lemma}
\label{span}
If $\langle R'\cup R''\rangle=R$ then $\dim(R')+\dim(R'')\geq \dim(R)-(1,1)$.
\end{lemma}

\begin{pf}
This is proved in a similar manner to Lemma \ref{travapps} (i).  If the sum of the widths of $R',R''$ is less than the width of $R$ minus one, then either $R$ has two adjacent columns which do not intersect $R'\cup R''$, or the East-most or West-most column of $R$ does not intersect $R'\cup R''$.  In either case, no site in such a column can be in $\langle R'\cup R''\rangle$.
\end{pf}

\begin{pfof}{lemma \ref{pod}}
The proof is by induction on the number of vertices of ${\cal H}$.  Suppose the lemma holds for all hierarchies with fewer vertices than ${\cal H}$.  We consider three cases according to whether the root $r$ is a seed, normal, or a splitter.

If $r$ is a seed (so it is the only vertex), then we take $S=R$, and the result holds trivially.

If $r$ is normal, so that $r\Rightarrow y$ say, then we apply the inductive hypothesis to the sub-hierarchy ${\cal H}'$ rooted at $y$ (that is, the hierarchy obtained by taking all vertices and edges in directed chains $y=v_0\leadsto v_1\leadsto v_2\leadsto \cdots $ away from $y$, together with the associated rectangles).  Let $S=S({\cal H})=S({\cal H}')$, and note that ${\cal H}'$ has the same number of splitters as ${\cal H}$, to obtain
\begin{eqnarray*}
\sum_{u\Rightarrow v} U(R_v,R_u) &\geq&
U(R_y, R) + U(S,R_y) - \nsplit 2qg(Z) \\
&\geq & U(S,R) - \nsplit 2qg(Z),
\end{eqnarray*}
by Proposition \ref{vsum} and the definition of $U(\cdot,\cdot)$.

If $r$ is a splitter, so that $r\rightrightarrows (y_1,y_2)$ say, we apply the inductive hypothesis to the sub-hierarchies ${\cal H}_1,{\cal H}_2$ rooted at $y_1,y_2$, and denote their pods $S_1=S({\cal H}_1), S_2=S({\cal H}_2)$.  We also write $R_1=R_{y_1}, R_2=R_{y_2}$.  Since the total number of splitters in ${\cal H}_1$ and ${\cal H}_2$ is one less than $\nsplit$, we obtain
\begin{equation}
\label{split1}
\sum_{u\Rightarrow v} U(R_v,R_u) \geq U(S_1,R_1)+U(S_2,R_2) -(\nsplit-1) 2qg(Z).
\end{equation}
Now we apply Proposition \ref{vsplit} with
\begin{eqnarray*}
\vc{a}&=&q \dim(S_1), \\
\vc{b}&=&q \dim(R_1), \\
\vc{c}&=&q \dim(S_2), \\
\vc{d}&=&q \dim(R_2), \\
\vc{r}&=&q \dim(R).
\end{eqnarray*}
We choose the pod $S=S({\cal H})$ to be a rectangle satisfying $S\subseteq R$ and
$$\vc{s}=q \dim(S),$$ 
where $\vc{s}$ is as in Proposition \ref{vsplit}. The formula for $\vc{s}$ ensures that the dimensions of $S$ are indeed integers and that $S$ depends only on ${\cal H}$.  Furthermore, since $\vc{s}\leq \vc{a}+\vc{c}$ we have
$$\dim(S)\leq\dim(S_1)+\dim(S_2)\leq \sum_{w \text{ seed}} \dim(R_w),$$
by the inductive hypothesis, since the set of seeds of ${\cal H}$ is the disjoint union of the sets of seeds of ${\cal H}'$ and ${\cal H}''$.
It is easy to check that the conditions of Proposition \ref{vsplit} are met, by Lemma \ref{span} and property (ii) of a good hierarchy, so we obtain
$$U(S_1,R_1)+U(S_2,R_2) \geq U(S,R) -2qg(Z).$$
Combining this with (\ref{split1}) gives the required bound.
\end{pfof}

Now we derive a lower bound on the second sum in (\ref{hocc}).
If $w$ is a seed then 
$$\frac{V(R_w)}{q\phi(R_w)}=\frac{\lng(R_w)\; g(q\sht(R_w))}{\lng(R_w)+\sht(R_w)}\geq\frac{g(2Z)}{2},$$
by property (i) of a good hierarchy.
Hence
\begin{equation}
\label{podbound}
\sum_{w \text{ seed}} V(R_w) \geq \frac{g(2Z)}{2}\sum_{w \text{ seed}} q\phi(R_w) \geq \frac{g(2Z)}{2} q\phi(S),
\end{equation}
since the pod $S$ satisfies $\dim(S)\leq \sum_{w \text{ seed}} \dim(R_w)$ by Lemma \ref{pod}.

Substituting from Lemma \ref{pod} and (\ref{podbound}) into (\ref{hocc}), for any good hierarchy ${\cal H}$ we have
\begin{eqnarray}
\label{hocc2}
\lefteqn{P_p({\cal H} \text{ occurs})\leq} \nonumber \\
&& Q^{\nnorm} Q_1^{\nsplit}
\exp-q^{-1}\left[(1-2c)U(S,R)+\frac{g(2Z)}{2} q\phi(S)\right],
\end{eqnarray}
where $Q_1=e^{2g(Z)}$.

Now suppose that ${\cal H}$ is a good hierarchy with root-label $R_r=R=R(\lceil B/q\rceil, \lceil B/q\rceil)$.  Let us find an upper bound on the number of vertices of ${\cal H}$.  By properties (iii),(v) of a good hierarchy, in any directed chain of vertices $r=v_0\leadsto v_1 \leadsto \cdots \leadsto v_k $ away from the root, at least half the edges have $\phi$ decreasing by at least $T/(2q)$, so the number of vertices in such a chain is at most
$$2 \frac{2 \lceil B/q\rceil}{T/(2q)}+1 \leq\frac{20B}{T}.$$
Hence, since the graph underlying ${\cal H}$ is a binary tree, the total number of vertices in ${\cal H}$ is at most
$$M:=2^{20B/T}.$$
All that matters is that this number depends only on $B$, not on $q$.

We now divide hierarchies into two different types according to the semi-perimeter of the pod.  If $q\phi(S)\leq A$ then by Lemma \ref{vexact} we have
$$U(S,R)\geq 2\int_{q\phi(S)}^{q\lceil B/q\rceil}g(z)\; dz \geq 2\int_A^B g(z)\; dz.$$
On the other hand, if $q\phi(S)>A$ then by the choice of $Z$, (\ref{zchoice}), we have that
$$\frac{g(2Z)}{2} q\phi(S) \geq \frac{g(2Z)A}{2} \geq 2 \lambda\geq 2(1-2c)\int_A^B g(z)\; dz,$$
since $\lambda=\int_0^\infty g(z) \; dz$.
In both cases we obtain from (\ref{hocc2}) that
\begin{equation}
\label{hocc3}
P_p({\cal H} \text{ occurs})\leq Q_2^M \exp -q^{-1}2(1-2c)\int_A^B g(z)\; dz,
\end{equation}
where $Q_2=\max\{Q,Q_1\}$.

We now bound the total number of possible good hierarchies with root-label $R$.  The number of abstract directed graphs with at most $M$ vertices is at most $M 2^{M^2}$, and the number of different rectangles in $R$ is at most $(B/q +1)^4$, so the number of possible good hierarchies is at most
\begin{equation}
\label{nhier}
M 2^{M^2}(2B/q)^{4M}.
\end{equation}

From Proposition \ref{hier} and (\ref{hocc3}), (\ref{nhier}) we deduce the existence of constants $K_1,K_2\in(0,\infty)$ depending only on $B$, such that when $q$ is sufficiently small,
$$P_p(R \text{ is internally spanned})\leq K_1 q^{-K_2} \exp -q^{-1}2\lambda_B,$$
where
\begin{equation}
\label{lambdab}
\lambda_B=(1-2c)\int_A^B g(z)\; dz.
\end{equation}

Hence, recalling that $q\sim p$, we have
$$\liminf_{p\rightarrow 0} -p\log I(\lceil B/q\rceil,p) \geq 2\lambda_B.$$
This implies that
\begin{equation}
\label{final}
\liminf_{p\rightarrow 0} -p\log I(\lfloor B/p\rfloor,p) \geq 2\lambda_B.
\end{equation}
To check this note that for any $\epsilon>0$ we may write $\lfloor B/p\rfloor=\lceil (B+\epsilon)/q'\rceil$, where $q'=-\log(1-p')$ and $p'\sim p(B+\epsilon)/B$ as $p\rightarrow 0$; since $I(L,p)$ is increasing in $p$ this implies that the left side of (\ref{final}) is at least $2\lambda_{B+\epsilon} B/(B+\epsilon)$, establishing (\ref{final}).  Finally, since $A=c=1/B$ we have from (\ref{lambdab}) that $$\lambda_B\rightarrow \int_0^\infty g(z)\; dz =\lambda \quad\text{ as } B\rightarrow \infty.$$
Hence we have proved Theorem \ref{mainprop} (ii). 

\section*{Acknowledgments}
I thank Marek Biskup, Lincoln Chayes, Nathaniel Grossman, Tom Liggett and Roberto Schonmann for valuable conversations.

\bibliography{phd}

\noindent
{\sc Alexander E. Holroyd}

\noindent
UCLA Department of Mathematics

\noindent
405 Hilgard Avenue

\noindent
Los Angeles

\noindent
CA 90095-1555

\noindent
U. S. A.

\noindent
{\tt holroyd@math.ucla.edu}

\end{document}